\documentclass[12pt]{amsart}

\usepackage{amsfonts}
\usepackage{amsmath}
\usepackage{amssymb}
\usepackage{mathabx}

\usepackage{mathrsfs}
\usepackage{amscd}
\usepackage[all]{xy}
\usepackage[pdftex,final]{graphicx}

\usepackage{hyperref}
\usepackage[hyperpageref]{backref}

\newtheorem{lem}{Lemma}[section]
\newtheorem{thm}[lem]{Theorem}
\newtheorem{prop}[lem]{Proposition}
\newtheorem{cor}[lem]{Corollary}
\theoremstyle{definition}
\newtheorem{defi}[lem]{Definition}

\newtheorem{exa}[lem]{Example}

\newtheorem{ques}[lem]{Question}
\newtheorem{rem}[lem]{Remark}

\renewcommand{\iff}{\Longleftrightarrow}

\newcommand{\Q}{\Bbb{Q}}
\newcommand{\F}[1]{\Bbb{F}_{#1}}

\newcommand{\Z}{\Bbb{Z}}
\newcommand{\N}{\Bbb{N}}
\newcommand{\R}{\Bbb{R}}

\newcommand{\K}{\mathcal{K}}

\newcommand{\LL}{L^{\tau}}

\renewcommand{\S}{\Sigma}

\newcommand{\wn}[1]{\mathcal{W}_{#1}}

\newcommand{\spl}[2]{\operatorname{SL}_{#1}\left(#2\right)}
\newcommand{\pspl}[2]{\mathrm{PSL}_{#1}(#2)}
\newcommand{\gl}[2]{\mathrm{GL}_{#1}(#2)}
\newcommand{\pgl}[2]{\mathrm{PGL}_{#1}(#2)}

\newcommand{\pb}[1]{\mathcal{P}(#1)}

\newcommand{\rpb}[1]{\mathcal{RP}(#1)}

\newcommand{\rpbker}[1]{\mathcal{RP}_1(#1)}

\newcommand{\qrpb}[1]{\widetilde{\mathcal{RP}}(#1)}
\newcommand{\hrpb}[1]{\widehat{\mathcal{RP}}(#1)}
\newcommand{\hpb}[1]{\widehat{\mathcal{P}}(#1)}

\newcommand{\crpb}[1]{\widecheck{\mathcal{RP}}(#1)}
\newcommand{\cpb}[1]{\widecheck{\mathcal{P}}(#1)}

\newcommand{\bl}[1]{\mathcal{B}(#1)}

\newcommand{\rbl}[1]{{\mathcal{RB}}(#1)}

\newcommand{\gpb}[1]{\left[ #1\right]} 

\newcommand{\suss}[2]{\psi_{#1}\left( #2\right)}
\newcommand{\Suss}[2]{\Psi_{#1}\left( #2\right)}
\newcommand{\Sus}[1]{\psi(#1)}

\newcommand{\ks}[2]{{\K}^{\mbox{\tiny $(#1)$}}_{#2}}

\newcommand{\bconst}[1]{C_{#1}}
\newcommand{\cconst}[1]{D_{#1}}

\newcommand{\spec}[1]{S_{#1}}

\newcommand{\pn}[2]{\mathbb{P}^{#1}(#2)}
\newcommand{\projl}[1]{\pn{1}{#1}}

\newcommand{\denom}[1]{\mathrm{denom}(#1)}

\newcommand{\card}[1]{\left| #1 \right|}

\renewcommand{\ker}[1]{\mathrm{Ker}(#1)}
\newcommand{\image}[1]{\mathrm{Im}(#1)}

\newcommand{\coker}[1]{\mathrm{Coker}(#1)}

\newcommand{\sqcl}[1]{G_{#1}}

\newcommand{\sgr}[1]{\mathrm{R}_{#1}}

\newcommand{\an}[1]{\left\langle{#1}\right\rangle}
\newcommand{\pf}[1]{\left\langle\!\left\langle{#1}\right\rangle\!\right\rangle}
\newcommand{\aug}[1]{\mathcal{I}_{#1}}

\newcommand{\igr}[1]{\Z[#1]}

\newcommand{\zhalf}[1]{{#1}\left[ \tfrac{1}{2}\right]}

\newcommand{\sgn}[1]{\mathrm{sgn}(#1)}

\newcommand{\kind}[1]{K^{\mathrm{\small ind}}_3(#1)}

\newcommand{\ho}[3]{\mathrm{H}_{#1}\left( #2,#3 \right)}

\renewcommand{\sgn}[1]{\mathrm{sgn}(#1)}

\title{Bloch Groups of local rings with residue field  $\F{3}$}
\author{  Kevin Hutchinson}
\address{School of Mathematics and Statistics,
 University College Dublin}
\email{  kevin.hutchinson@ucd.ie}
\date{\today}

\keywords{
$K$-theory, Group Homology
}
\subjclass{19F99, 20G10}

\begin{document}
\maketitle
\tableofcontents

\begin{abstract} We give a generic presentation of the pre-Bloch group of a local ring with residue field $\F{3}$. We apply this to the calculation of pre-Bloch groups and Bloch groups of some local rings.
\end{abstract}
\section{Introduction}\label{sec:intro}

The scissors congruence group, or pre-Bloch group, $\pb{F}$ of a field $F$, and its subgroup  the Bloch group, $\bl{F}$, (see \cite{sah:dupont}) is of interest because it describes the indecomposable $K_3$ of the field, $\kind{F}$,  modulo some known torsion coming from the group $\mu_F$ of roots of unity (\cite{sus:bloch}). 
 The refined scissors congruence group $\rpb{F}$  and the refined Bloch group $\rbl{F}$ were introduced in \cite{hut:cplx13} to help to understand the kernel of the natural 
surjective homomorphism $\ho{3}{\spl{2}{F}}{\Z}\to\kind{F}$. These groups and their properties were then used to calculate $\ho{3}{\spl{2}{F}}{\Z[\frac{1}{2}]}$ for local fields $F$ as well as well as for $F=\Q$ (\cite{hut:rbl11},\cite{hut:sl2Q}). The definitions and techniques in these results extend readily from fields to local rings with sufficiently large residue fields and thus allow us to describe $\kind{A}$ and $\ho{3}{\spl{2}{A}}{\Z[\frac{1}{2}]}$ for many local rings $A$ (\cite{mirzaii:bwlocal},\cite{hut:slr}, \cite{hut:mokarimirzaii}, \cite{mirzaii:perez25}) . 
However, the given funtorial definition of the (refined) pre- Bloch group does not extend to the fields $\F{2}$ and $\F{3}$ or to local rings with small residue fields.

In \cite{hut:bloch1} the author and  Rodrigo Cuitun Coronado defined functorially the pre-Bloch groups $\rpb{A}$ (and $\pb{A}$)  and the Bloch groups $\rbl{A}$ (and $\bl{A}$) of \emph{any} commutative ring $A$.   Under suitable hypotheses on the ring $A$, there is a natural homomorphism\\
 $\ho{3}{\spl{2}{A}}{\Z}\to \rpb{A}$ whose image is $\rbl{A}$ (or from $\ho{3}{\gl{2}{A}}{\Z}$ to $\pb{A}$ with image $\bl{A}$).  Those suitable hypotheses, which are satisfied by all fields and local rings, concern the low-dimensional homology of a simplicial complex $Y(A)$ associated to the ring $A$, introduced in \cite{hut:ge2}, and whose fundamental group can be described in terms of the algebraic $K$-theory of $A$ (\cite[Theorem 6.9]{hut:ge2}).

When $A$ is a field with at  least four elements or, more generally a local ring whose residue field has at least four elements, then $\pb{A}$ has a natural presentation in terms of generators $\gpb{u}$, $u\in \wn{A}:=\{ a\in A^\times|\  1-a\in A^\times\}$ (the set of exceptional units), and a family of $5$-term dilogarithm relations  $R_{u,v}$ where $u,v\in\wn{A}$ also satisfy $u/v\in \wn{A}$. Likewise, $\rpb{A}$ has an analogous presentation as an $\sgr{A}:=\Z[A^\times/(A^\times)^2]$-module. These presentations have proved useful in calculating $K$-theory groups and the third homology of $\spl{2}{A}$ (see, in particular, \cite{hut:sl2Q}).

 When $A$ is a ring mapping to $\F{2}=\Z/2$ the set $\wn{A}$ is necesssarily empty and the proposed generators $\gpb{u}$ do not exist. In the (prior) companion paper to the current paper (\cite{hut:f2}), we showed how to find a complete set of generators and relations of $\rpb{A}$ when $A$ is a local ring with residue field $\F{2}$. 

When $A$ admits a homomorphism to $\F{3}=\Z/3$ and if $u,v\in\wn{A}$ then $u/v\not\in \wn{A}$. Thus when $A$ is a local ring with residue field $\F{3}$, one can show that $\rpb{A}$ is generated  by the set $\{ \gpb{u}\ |\ u\in\wn{A}\}$ but that the above $5$-term relations no longer exist. In this paper, we show how to write down a complete set of relations among the generators $\gpb{u}$, $u\in\wn{A}$, of $\rpb{A}$ when $A$ is a local ring with residue field $\F{3}$ (Proposition \ref{prop:rpbapresf3}).

In section \ref{sec:review} we review some needed results about simplicial complexes in general and the about simplicial complex $Y(A)$, whose low-dimensional homology plays a key role.

If $A$ is a local ring with residue field $\F{3}$ then $Y(A)$ is a $3$-dimensional simplicial complex (whose geometric realization is homotopy equivalent to a pointed sum of $3$-spheres, by Proposition \ref{prop:localhom}). In section \ref{sec:3dim}, we establish some results about groups of cycles the complex $L_\bullet$ of non-degenerate chains associated a $3$-dimensional simplicial complex $Y$.

In section \ref{sec:hrpba} we introduce the module $\hrpb{A}$, which maps naturally to $\rpb{A}$ for any ring $A$. It arises naturally from the algebra of the complex $L_\bullet(A)$ studied in section \ref{sec:3dim}.   We show that for a large class of rings --  including local rings with residue field $\F{3}$  -- it has a presenation with generators $\gpb{u}$, $u\in \wn{A}$ and relations consisting of known generic identities involving the special classes of elements $\Suss{i}{u}$ and $\cconst{A}(u)$, already enumerated in \cite{hut:bloch1}. For local rings with residue field $\F{3}$,  there is an exact sequence \\
$H_3(|Y(A)|)\to\hrpb{A}\to\rpb{A}\to 0$. We calculate the modules $\rpb{A}$ when $A=\Z/3$ and $A=\Z/9$.

In section \ref{sec:presrpa} we give a presentation of the module $\rpb{A}$) -- and the abelian group $\pb{A}$ -- when $A$ is a local ring with residue field $\F{3}$. Furthermore, we describe the structure of the cokernel of the natural map $\crpb{A}\to\hrpb{A}$ where $\crpb{A}$ is the module introduced in \cite{hut:f2}.

In the final section, \ref{sec:exa}, we consider some examples. In particular, we investigate the injectivity or otherwise of the map $\hrpb{A}\to\rpb{A}$.  That is, in the absence of the dilogarithm identities, what additional relations among the generators $\gpb{u}$ are needed, beyond the generic algebraic identities of \cite{hut:bloch1}, to define $\rpb{A}$. We show that the natural map $\hrpb{\Z[\frac12]}\to \rpb{\Z[\frac12]}$ is an isomorphism. We show that the kernel of the map $\hrpb{\Z_{\an{3}}}\to \rpb{\Z_{\an{3}}}$ contains a certain free abelian group of infinite rank with explicit basis. 

\subsection{Notation}  In this article all  rings are commutative and have a unit.

For a ring $A$, $A^\times$ denotes the group of units. $\sgr{A}$ denotes the group ring
$\Z[A^\times/(A^\times)^2]$ of the group of square classes of units. For a unit $u\in A^\times$, the image in $\sgr{A}$ will be denoted $\an{u}$. The augmentation ideal $\aug{A}$ of the ring $\sgr{A}$ is the additive subgroup generated by the elements $\pf{u}:=\an{u}-1$, $u\in A^\times$.  

\section{Review}\label{sec:review}
In this section we review some terminology and needed results, primarily from \cite{hut:f2}.
\subsection{Homology of Simplicial Complexes}
Let $Y$ be a simplicial complex.  For $n\geq 0$, $Y_n$ will denote the set of $n$-simplices of $Y$. Let $C_\bullet=C_\bullet(Y)$  be the complex whose $n$-th term $C_n$ is the free abelian group on the set 
\[
E_n=E_n(Y):=\{ (x_0,\ldots,x_n)\ |\ \{ x_0,\ldots,x_n\} \mbox{ is an $m$-simplex  for some }m\leq n\},
\] 
equipped with the usual boundary homomorphism 
\[
d_n (x_0,\ldots,x_n):= \sum_{i=0}^n(-1)^i(x_0,\ldots, \widehat{x_i},\cdots,x_n).
\]
Thus $H_n(C_\bullet)\cong H_n(|Y|)$ for all $n\geq 0$, where $|Y|$ is the geometric realization of $Y$.
Let $L_\bullet$ be the subcomplex whose $n$-the term is $\Z[X_n]$ where $X_n$ is the set of (non-degenerate) ordered $n$-simplices:
\[
X_n=X_n(Y):=\{ (x_0,\ldots,x_n)\ |\ \{ x_0,\ldots,x_n\} \mbox{ is an $n$-simplex}\}.
\]
Let $\bar{C}_\bullet$ denote the quotient complex $C_\bullet/L_\bullet$. Thus, as an abelian group, $\bar{C}_n$ is isomorphic to  the free abelian group $D_n=D_n(Y)$ on the  set 
$E_n\setminus X_n$ of (ordered) degenerate simplices and $C_n=L_n\oplus D_n \cong L_n\oplus \bar{C}_n$ as abelian groups for all $n$.

If $Y$ is a simplicial complex and if $\alpha\in Y_m$ we will denote by $Y^\alpha$ the subcomplex defined by $Y^\alpha_n:=\{ \beta\in Y_n\ |\ \alpha\cup\beta\in Y_{n+m+1}\}$. A path or loop in $Y^\alpha$ will also be called an $\alpha$-path or $\alpha$-loop in $Y$. Thus $a_0,a_1,\ldots,a_n$  ($a_i\in Y_0$) is an $\alpha$-path in $Y$ if $\{ \alpha_{i-1},\alpha_i\}\cup\alpha\in Y_{m+2}$ for $i=1,\ldots,n$.  

For $m\geq 0$, we will say that $Y$ \emph{satisfies (the condition) $(C)_m$}  if $Y^\alpha$ is path connected for all $\alpha\in Y_m$.

Below we will use:
\begin{lem}(\cite[Lemma 3.2]{hut:f2})\label{lem:c2} Suppose that every $1$-simplex of  $Y$  is contained in a $2$-simplex. Then $H_2(\bar{C}_\bullet)=0$.
\end{lem}
\subsection{The simplicial complex $Y(A)$}
Let $A$ be a commutative ring.  Let $Y_0(A)$ be the set of unimodular rows $u\in A^2$ up to multiplication by an element of $A^\times$. Furthermore, for any $a\in A$, we will also denote the unimodular row $(a,1)\in Y_0(A)$ by $a$. For any $u\in A^\times$, the symbol $\infty$ will denote the class of $(u,0)$ in $Y_0(A)$. We also set $\infty:=(1,0)\in  Y_0(A)$. Let $Y:=Y_n(A)$ be the simplicial complex whose set $Y_n(A)$ of $n$-simplices consists of all sets $\{ [u_0],\ldots,[u_n]\}$, $[u_i]\in Y_0(A)$ with the property that the matrices
\[
M_{ij}:=
\left[
\begin{array}{c}
u_i\\
u_j\\
\end{array}
\right]
\]
lie in $\gl{2}{A}$. The group $\pgl{2}{A}$, and therefore also $\gl{2}{A}$, acts on $Y(A)$ by right multiplication. The associated complex $L_\bullet=L_\bullet(Y(A))$ is naturally a complex of $\pgl{2}{A}$-modules.

 We are concerned with determining $Z_2(L_\bullet)$ using the results of section (\ref{sec:intro}).
More particularly we would like to calculate the groups   $\rpb{A}:=(Z_2(L_\bullet))_{\pspl{2}{A}}$ and $\pb{A}:=(Z_2(L_\bullet))_{\pgl{2}{A}}$.
 Because of the exact sequence
 \[
1\to\pspl{2}{A}\to\pgl{2}{A}\to A^\times/(A^\times)^2\to 1,
\]
 the group $\rpb{A}$ is naturally an $\sgr{A}$-module (where $\sgr{A}:= \Z[A^\times/(A^\times)^2]$): for all $z\in Z_2(L_\bullet)$,  $\an{a}\cdot [z] = [z\cdot M]$ for any $M\in \gl{2}{A}$ with determinant $a$, and where 
$[z]$ denotes the class of $z$ in  $\rpb{A}$.  When $M\in\pgl{2}{A}$, $\an{\mathrm{det}(M)}$ will denote the resulting square class in $A^\times/(A^\times)^2$.  If $N$ is a (right) 
$\pgl{2}{A}$-module, then the group of coinvariants $M_{\pspl{2}{A}}$ is thus naturally an $\sgr{A}$-module via
\[
\an{\mathrm{det}(M)}\cdot\bar{x}=\overline{x\cdot M}
\]
for all $x\in N$, $M\in \pgl{2}{A}$.

Note  that $\{ 0,\infty,1\}\in Y_2(A)$ for any ring $A$. Virtually by definition, the group $\pgl{2}{A}$ acts transitively on $X_1(A)$, and hence on the set of $1$-simplices $Y_1(A)$; i.e., if $a=(a_1,a_2), b=(b_1,b_2)$ are the rows of an invertible matrix, then the resulting matrix 
\[
T_{a,b}:=\left[
\begin{array}{cc}
a_1&a_2\\
b_1&b_2\\
\end{array}
\right]\in \pgl{2}{A}
\]
has the property that $(\infty,0)\cdot T_{a,b}=(a,b)$.
 Furthemore, $\gl{2}{A}$ acts transitively on $X_2(A)$, and hence also on $Y_2(A)$:
Let $\{  a=(a_1,a_2),b=(b_1,b_2),c=(c_1,c_2)\}\in X_2(A)$. Let
\[
M_{a,b,c}:= 
\left[
\begin{array}{cc}
d(b,c)&0\\
0&d(c,a)\\
\end{array}
\right]\cdot T_{a,b}=
\left[
\begin{array}{cc}
a_1d(b,c)&a_2d(b,c)\\
b_1d(c,a)&b_2d(c,a)\\
\end{array}
\right]
\in \pgl{2}{A}
\]
where, whenever $x=(x_1,x_2), y=(y_1,y_2)$ are the rows of a $2\times 2$ invertible matrix, we define
\[
d(x,y):= \mathrm{det}
\left[
\begin{array}{cc}
x_1&x_2\\
y_1&y_2\\
\end{array}
\right].
\]
Then $(\infty,0,1)\cdot M_{a,b,c}=(a,b,c)$. 

Observe that for any $(a,b,c)\in X_2(A)$, we have $\an{\mathrm{det}(M_{a,b,c})}=\an{d(a,b)d(b,c)d(c,a)}\in A^\times/(A^\times)^2$.

We can express properties of the  simplicial complex $Y(A)$ in direct ring-theoretic terms. For a ring $A$ we let $\wn{A}:=\{ a\in A\ |\ a(1-a)\in A^\times\}\subset A^\times$. 
\begin{lem}\label{lem:ya} Let $A$ be a ring.
\begin{enumerate}
\item Every $0$ -simplex of $Y(A)$ is contained in a $1$-simplex and every $1$-simplex is contained in a $2$-simplex.
\item $Y_3(A)\not=\emptyset$ if and only if $\wn{A}\not=\emptyset$ and, if so, every $2$-simplex is contained in a $3$-simplex.
\item The simplicial complex $Y(A)$ satisfies the condition $(C)_1$ if and only if $\wn{A}$ generates $A^\times$ as a group.
\end{enumerate}
\end{lem}

\begin{proof} 
\begin{enumerate}
\item Every $0$-simplex is contained in a $1$-simplex by the definition of \emph{unimodular}. If $\{a,b\}$ is a $1$-simplex, then it is contained in the $2$-simplex $\{ a,b,c\}$ where 
$c=(1 1)\cdot T_{a,b}^{-1}$. 
\item Note that $Y^{\{ \infty,0\}}_0=\{ (u,1)\ |\ u\in A^\times\}=A^\times$. Now one easily verifies that if $u\in A^\times$, then $\{ \infty,0,1,u\}\in Y_3$ if and only if $1-u\in A^\times$; i.e., 
if and only if $u\in \wn{A}$. Conversely, if $\{ a,b,c,d\}$ is any $3$-simplex, then so is $\{\infty,0,1,u\}$ where $u=d\cdot M_{a,b,c}^{-1}$. 

Conversely,  if $u\in \wn{A}$ and if $\{ a,b,c\}\in Y_2(A)$, then $\{a,b,c,d\}\in Y_3(A)$ when $d=u\cdot M_{a,b,c}$. 
\item Since $\gl{2}{A}$ acts transitively on edges of $Y(A)$, it follows that for any edge $\{ x,y\}\in Y_1(A)$, the simplicial complex $Y^{\{ x,y\}}$ is path-connected if and only if 
$Y^{\{ \infty,0\}}$ is path-connected. Recall that $Y^{\{ \infty,0\}}_0=A^\times$. Furthermore, if $u,v\in Y^{\{ \infty,0\}}_0$, then $\{ u,v\}\in Y^{\{ \infty,0\}}_1$ if and only if 
$\{ 1,u^{-1}v\}\in Y^{\{ \infty,0\}}_1$ (since multiplication by $\mathrm{diag}(u^{-1},1)$ is an automorphism of $Y^{\{ \infty,0\}}$) if and only if $\{ \infty,0,1,u^{-1}v\}\in Y_3(A)$ if and only if $u^{-1}v\in \wn{A}$. Thus an $\{ \infty,0\}$-path from $u$ to $v$ in $Y(A)$ is a sequence of units $u=u_0,u_1,\cdots, u_{n-1},u_n=v$ satisfying $w_i:=u_{i-1}^{-1}u_i\in \wn{A}$ for 
$i=1,\ldots,n$. It follows that $u^{-1}v= w_1\cdots w_n$. 

Conversely, given $u,v\in A^\times$, if there exist $w_1,\ldots,w_n\in\wn{A}$ such that $u^{-1}v=w_1w_2\cdots w_n$, then $u, uw_1,uw_1w_2,\cdots,uw_1\cdots w_n=v$ is an $\{ \infty,0\}$-path 
in $Y(A)$ from $u$ to $v$. 
\end{enumerate}
\end{proof}

\begin{rem} Similarly, it can also be easily seen that $Y(A)$ has property $(C)_0$ if and only if $A^\times$ generates $A$ as an additive group, but we will not need this property 
below.
\end{rem}

\subsection{The connection with the low-dimensional homology of $\spl{2}{A}$}
Let $\LL_\bullet:=\LL_\bullet(A)$ be the following truncation of the complex $L_\bullet$:
\[
\LL_n(A)=\left\{
\begin{array}{ll}
L_n(A),& n\leq 2\\
Z_2(L_\bullet):=\ker{L_2(A)\to L_1(A)},& n=3\\
0,& n>3
\end{array}
\right.
\]
(where the map $\LL_3(A)=Z_2(L_\bullet)\to \LL_2(A)=L_2(A)$ is the inclusion homomorphism).
 Let $\epsilon:L_0(A)=Z[X_0]\to \Z$ be the homomorphism sending each $x\in  X_0$ to $1$. If we regard $\Z$ as a complex of (trivial) $\pgl{2}{A}$-modules concentrated in dimension $0$, then $\epsilon$ induces maps of complexes $\epsilon:\LL_\bullet \to \Z$ and $\epsilon:L_\bullet\to \Z$. It can be shown (\cite{hut:ge2}) that this map induces an isomorphism on homology for many rings $A$,
including when $A$ is a field or a local ring or when $A=\Z$. It follows that for such rings $A$, and for any group $G$ admitting a homomorphism to $\pgl{2}{A}$, there is a hyperhomology spectral sequence of the form 
\[
E^1_{p,q}=\ho{p}{G}{L_q}\Longrightarrow \ho{p+q}{G}{\LL_\bullet}\cong\ho{p+q}{G}{\Z}.
\]

Thus, when $G=\spl{2}{A}$, we have $E^1_{0,3}=\ho{0}{\spl{2}{A}}{Z_2(L_\bullet)}=\rpb{A}$ and $E^1_{0,2}=\ho{0}{\spl{2}{A}}{L_2}\cong\sgr{A}$ via $(a,b,c)\mapsto \an{d(a,b,c)}$.

 We let $\lambda_1$ be the map 
\[
d^1_{0,3}:\rpb{A}=E^1_{0,3}\to E^1_{0,2}=\sgr{A}
\]
and we let $\rpbker{A}:=\ker{\lambda_1}=E^2_{0,3}$. Now $E^2_{1,1}=0$ always (\cite[Corollary 4.4]{hut:bloch1}). We let $\lambda_2$ be the map $d^3_{0,3}:\rpbker{A}\to E^3_{2,0}$. Then $E^4_{0,3}=E^\infty_{0,3}=\ker{\lambda_2:\rpbker{A}\to E^3_{2,0}}=\rbl{A}$ is the \emph{(refined) Bloch group} of the ring $A$. Thus when $A$ is a field or a local ring or $\Z$, there is a natural surjective edge homomorphism $\ho{3}{\spl{2}{A}}{\Z}\to \rbl{A}\subset\rpb{A}$. (For more details about this spectral sequence, see \cite[Sections 4,5]{hut:bloch1}
\subsection{The space $Y(A)$ for local rings $A$}
For a ring $A$ let $\Gamma(A)$ be the $1$-skeleton of the simplicial complex $Y(A)$. Thus $\Gamma(A)$ is a graph whose vertices are classes of unimodular rows and whose edges are pairs 
$\{ [u_0], [u_1]\}$ with the property that the matrix\\
 $T_{u_0,u_1}=\left[
\begin{array}{c}
u_0\\
u_1\\
\end{array}
\right]$ is invertible. Then, by definition, the simplicial complex $Y(A)$ is just the clique complex of the graph $\Gamma(A)$; i.e., it is the complex whose $r$-simplices are the vertices of $(r+1)$-cliques (complete subgraphs on $r+1$ vertices)  of $\Gamma(A)$.

Let $A$ be a local ring with maximal ideal ${M}$ and residue field $k=A/M$. Then $(a_1,a_2)$ is unimodular if and only if either $a_1\in A^\times$ or $a_2\in A^\times$. If $a_2\in A^\times$, 
then $(a_1,a_2)=(a,1)=a\in A$ where $a=a_1a_2^{-1}$. If $a_1\in A^\times$, then $(a_1,a_2)=(1,b)$ where $b=a_1^{-1}a_2\in A$. For convenience, we will use the notation $1/b$ for such classes $(1,b)$ of unimodular rows.  Note that $\infty:=1/0$ in this notation. Also,  when $a\in A^\times$ we have $1/a=(1,a)=(a^{-1},1)=a^{-1}$. It follows that $1/A\cap A=A^\times$ and the vertices of $\Gamma(A)$ 
are $A\cup 1/A=A\sqcup 1/M$. 

For $a,b\in A\subset \Gamma(A)_0$, $\mathrm{det}(T_{a,b})=a-b$ and hence $\{ a,b\}$ is an edge if and only if $a-b\in A^\times$ if and only if $a\not\equiv b\pmod{M}$.
For any $a\in A, m\in M$, $\mathrm{det}(T_{a,1/m})=am-1\in A^\times$ and hence $\{ a,1/m\}$ is always an edge. Finally, for 
any $m,n\in M$, $\mathrm{det}(T_{1/m,1/n})=n-m\in M$ and hence $\{ 1/m,1/n\}$ is never an edge of $\Gamma(A)$.

Suppose now that the residue field $k$ is finite of cardinality $q$.
For $i\in \projl{k}$ define $T_i\subset Y(A)_0$ by
\[
T_i=\left\{
\begin{array}{ll}
\mbox{ the corresponding coset of $M$ in $A$},& i\in k\\
1/M,& i=\infty.\\
\end{array}
\right.
\]  
Note that each $T_i$ is naturally in bijective correspondence with $M$ for all $i\in\projl{k}$. In particular, if $A$ is finite, then $|T_i|=|M|=|A|/q$ for all $i\in\projl{k}$.

We can describe the
 vertices, $\Gamma(A))_0$, and edges, $\Gamma(A)_1$, of the graph $\Gamma(A)$ as follows:\\
 $\Gamma(A)_0=\bigsqcup_{i\in\projl{k}}T_i$ and 
$\{ t,s\}\in \Gamma(A)_1$ if and only if $t\in T_{i}$ and $s\in T_{j}$ with $i\not=j$; in other words, $\Gamma(A)$ is the \emph{complete multipartite 
graph}, $K(\{ T_i\}_{i\in\projl{k}})$  on the collection of sets $\{ T_i\}_{i}$. It follows that $Y(A)=Y(\{ T_i\}_{i})$ is the clique complex of this complete multipartite 
graph. In detail, $\{ t_0,\ldots,t_r\}\subset \Gamma(A)_0$ is an $r$-simplex of $Y(A)$ if and only if $t_k\in T_{i_k}$, for $0\leq k\leq r$ with $i_k\not=i_{k'}$ when $k\not=k'$. 
To review:

\begin{prop}\label{prop:localya} Let $A$ be a local ring with maximal ideal $M$ and finite residue field $K$ of size $q$. Then $\Gamma(A)$ is the complete multipartite graph 
$K(\{ T_i\}_{i\in\projl{k}})$ where the $T_i$
for $i\in k$ is the corresponding  coset of $M$  in $A$ and $T_\infty:= 1/M\subset Y(A)_0$. Thus the simplicial complex $Y(A)$ is the clique complex of this graph.
\end{prop}


For example, for any prime $p$, $\Gamma(\Z/p^n)$ is the complete bipartite graph $K(\{T_i\}_{i=0}^p)$ on $p+1$ sets $T_i$, each of cardinality $|M|=p^{n-1}$; it is the complete 
multipartite graph 
\[
\Gamma(\Z/p^n)\cong K_{\underbrace{p^{n-1},\ldots,p^{n-1}}_{p+1}}
\]
and hence $Y(\Z/p^n)$ is the corresponding clique complex
\[
Y(\Z/p^n)\cong Y_{\underbrace{p^{n-1},\ldots,p^{n-1}}_{p+1}}.
\]

\begin{prop}\cite[Proposition 2.5]{hut:f2}\label{prop:localhom} Let $A$ be a local ring with maximal ideal $M$ and a finite residue field $k=A/M$ of size $q$. 
For convenience, fix a bijection between $\projl{k}$ and $[0,q]:=\{ i\in \Z\ | \ 0\leq i\leq q\}$. For $i\in [0,q]$, let $T_i$  be as above.
For each $i\in [0,q]$, fix $x_i\in T_i$.
\begin{enumerate} 
\item For any $s=(s_0,\ldots,s_q)\in \prod_{i=0}^q (T_i\setminus\{ x_i\})$ let $Y_s$ be the simplicial complex\\
 $Y(\{ x_0,s_0\};\cdots ;\{ x_q,s_q\})$; i.e., the clique complex of the complete multipartite graph on this collection of  $q+1$ sets. Then $|Y_s|$ is homeomorphic to a $q$-sphere and the cycle
 \[
\gamma^s_q:=\sum_{y=(y_0,\ldots,y_q)\in\prod_{i=0}^q\{ x_i,s_i\}}(-1)^{\delta(y,s)}y\in  L_q(Y_s)\subset C_q(Y_s)
\]
(where $\delta(y,s):=\card{\{ i\ |\ y_i\not=s_i\}}=\card{\{ i\ |\ y_i=x_i\}}$) represents a generator of $H_q(|Y_s|)\cong \Z$.
\item The inclusion maps of simplicial complexes $Y_s\to Y(A)$, as $s$ varies over\\
 $S:=\prod_{i=0}^q (T_i\setminus\{ x_i\})$,  induce a map of simplicial complexes $\bigvee_s Y_s\to Y(A)$ 
(where the $Y_s$ are pointed with basepoint $x_0$) which induces
a homotopy equivalence of spaces
\[
\bigvee_{s\in S}|Y_s|\cong |Y(A)|.
\]
In particular, $H_n(|Y(A)|)=0$ for $n\not=0,q$ and $H_q(|Y(A)|)\cong \oplus_{s\in S} \Z$ with basis $\{ \gamma_q^s\ |\ s\in S=\prod_{i=0}^q (T_i\setminus\{ x_i\})\}$.
\end{enumerate}
\end{prop}
 Thus in this case $Y(A)$ is homotopy equivalent to a pointed sum of $q$-spheres and its $q$-dimensional homology is a free $\Z$-module with an explicit basis. For example, by Proposition \ref{prop:localhom}, $Y(\Z/p^n)$ has the homotopy type of a 
pointed sum of $(p^{n-1}-1)^{p+1}$ $p$-spheres. 
\subsection{Review of some  results in low-dimensions}
For $\sigma\in S_{n+1}$ and $\{ x_1,\ldots,x_{n+1}\}\in Y_n$ we denote by $R_{\sigma}(x_1,\ldots,x_{n+1})$ the element $(x_1,\ldots,x_{n+1})-\sgn{\sigma}(x_{\sigma(1)},\ldots,x_{\sigma(n+1)})$ 
in $L_n$.  For $n\geq 0$, let $\Sigma_n\subset L_n$ be the subgroup generated by the elements $R_{\sigma}(x_1,\ldots,x_{n+1})$, $\{x_1,\ldots,x_{n+1}\}\in Y_n$, $\sigma\in \S_{n+1}$.

\begin{lem}\cite[Lemma 3.6]{hut:f2} \label{lem:gen}  For any $n\geq 1$, $\Sigma_n$ is generated by elements of the form
 \[
R_{(i\  i+1)}(x_1,\ldots,x_{n+1})=(x_1,\ldots, x_i,x_{i+1},\ldots ,x_{n+1})+(x_1,\ldots, x_{i+1},x_i\ldots ,x_{n+1})
\]
with $1\leq i\leq n$, $\{ x_1,\ldots,x_{n+1}\}\in Y_n$.
\end{lem}

\begin{proof}  This follows from the identity
\[
R_{\sigma\tau}(x_1,\ldots,x_{n+1})=R_{\tau}(x_1,\ldots,x_{n+1})+\sgn{\tau}R_{\sigma}(x_{\tau(1)},\ldots,x_{\tau(n+1)})
\]
together with the fact that the symmetric group  $S_{n+1}$ is generated by the set\\
 $\{ (i\ i+1)\ |\ 1\leq i\leq n\}$.
\end{proof}

If  $\{ x,y\}\in Y_1$ and if $\{ x,y\}\cup \{ a_1,\ldots,a_{n-1}\}\in Y_n$  we will let 
 \begin{eqnarray*}
R^{\{ x,y\}}_i(a_1,\ldots,a_{n-1}):&=& R_{(i\ i+1)}(a_1,\ldots,a_{i-1},x,y,a_i,\ldots,a_{n-1})\\
:&=& (a_1,\ldots,a_{i-1},x,y,a_i,\ldots,a_{n-1})+(a_1,\ldots,a_{i-1},y,x,a_i,\ldots,a_{n-1}).\\
\end{eqnarray*}
Thus $\Sigma_n$ is generated by the elements of the form $R^{\{ x,y\}}_i(a_1,\ldots,a_{n-1})$.

The following lemma can be verified directly from the definitions:
\begin{lem}\label{lem:d} If $\{x,y\}\cup\{ a_1,\ldots,a_{n-1}\}\in Y_n$ then for all $1\leq i\leq n$ we have
\[
d_n(R^{\{ x,y\}}_i(a_1,\ldots,a_{n-1}))=\sum_{j=1}^{i-1}(-1)^{j-1}R^{\{x,y\}}_{i-1}(a_1,\ldots,\hat{a_j}\ldots,a_{n-1})+\sum_{j=i}^{n-1}(-1)^{j-1}R^{\{x,y\}}_{i}(a_1,\ldots,\hat{a_j}\ldots,a_{n-1})
\]
(where the first term of the right is the empty sum  - i.e., equals $0$ --  if $i=1$, and the second term is $0$ if $i=n$).
\end{lem}

In view of Lemmas \ref{lem:gen} and \ref{lem:d} we immediately deduce:

\begin{cor}\label{cor:subcplx}  For all $n\geq 1$, $d_n(\Sigma_n)\subset \Sigma_{n-1}$ and hence $\Sigma_\bullet:=\{ \Sigma_n,d_n\}_{n\geq 0}$ is a subcomplex of $L_\bullet$.
\end{cor}

We may thus set $Z_n(\Sigma_\bullet):= \ker{d_n:\Sigma_n\to \Sigma_{n-1}}$ for $n\geq 1$. Of course, $Z_n(\Sigma_\bullet)\subset Z_n=Z_n(L_\bullet)$ for all $n\geq 1$.

\begin{prop}\cite[Proposition 3.9]{hut:f2}\label{prop:symm} 
 For all $n\geq 1$, $d_{n+1}(D_{n+1})\cap L_n\subset Z_n(\Sigma_\bullet)$.
\end{prop}


We introduce some families of elements in $Z_2(\Sigma_\bullet)$:  First observe that whenever $\{ x,y,a\}\in Y_2$ we have \[
d_2(R_1^{\{ x,y\}}(a))=d_2(R_2^{\{ x,y\}}(a))=R_1^{\{ x,y\}}:= (x,y)+(y,x).
\]
Whenever $\{x,y,a\},\{x,y,b\}\in Y_2$, we have the following elements of  $Z_2(\Sigma_\bullet)$:
\begin{eqnarray*}
D_1^{\{ x,y\}}(a,b)&:=& R_1^{\{ x,y\}}(b)-R_1^{\{ x,y\}}(a)\\
D_2^{\{ x,y\}}(a,b)&:=& R_2^{\{ x,y\}}(b)-R_2^{\{ x,y\}}(a).\\
\end{eqnarray*}
Furthermore, for any $\{ x,y,a\}\in Y_2$, let $C^{\{ x,y\}}(a):=R_1^{\{ x,y\}}(a)-R_2^{\{ x,y\}}(a)\in Z_2(\Sigma_\bullet)$. 






Now $\Sigma_2=\Sigma_2^1+\Sigma_2^2$ where $\Sigma_2^i$ is the group generated by the terms $R_i^{\{ x,y\}}(a)$ for $\{ x,y,a\}\in Y_2$ and $i=1,2$.

For $i=1,2$, let $Z^i\subset \Sigma_2^i$ be the subgroup generated by the terms $D_i^{\{ x,y\}}(a,b)$ and let $Z^c$ be the subgroup generated by the terms $C^{\{ x,y\}}(a)$.

\begin{prop}\cite[Corollary 3.17]{hut:f2}\label{cor:z2pres}
$Z_2(\Sigma_\bullet)=Z^c\oplus Z^1=Z^c\oplus Z^2$  where $Z^c$ has the presentation:

Generators $C^{\{ x,y\}}(a)$ , $\{ x,y,a\}\in Y_2$\\
Relations:  For all $\{ x,y,z\}\in Y_2$, $C^{\{ x,y\}}(z)+C^{\{ x,z\}}(y)+C^{\{ y,z\}}(x)=0$, and

$Z^i$ has the presentation

Generators: $D_i^{\{ x,y\}}(a,b)$   with $\{ x,y\}\in Y_1$, $a,b\in Y^{\{ x,y\}}_0$\\
 Relations: $D_i^{\{ x,y\}}(a,b)+D_i^{\{ x,y\}}(b,c)=D_i^{\{ x,y\}}(a,c)$\\

Here we allow $a=b$, $b=c$ or $a=c$ and define $D_1^{\{ x,y\}}(z,z):=0$.
\end{prop}

\subsection{$Z_2(\Sigma_\bullet(A))$ and $\crpb{A}$}
We specialize to the case  $Y=Y(A)$ where $A$ is a ring. Then $Z_2(\Sigma_\bullet(A))$  is a $\pgl{2}{A}$-module and we define $\crpb{A}:=Z_2(\Sigma_\bullet(A))_{\spl{2}{A}}$.
Thus there is a natural homomorphism $\crpb{A}\to \rpb{A}=Z_2(L_\bullet(A))_{\spl{2}{A}}$ which is studied in \cite{hut:f2}.

We have the special elements $\suss{i}{u}:=D_i^{\{ \infty, 0\}}(1,u)$, $u\in A^\times$ and $\cconst{A}:=C^{\{ \infty,0\}}(1)$ in $Z_2(\Sigma_\bullet(A))$.  
Let $\omega:=
\left[
\begin{array}{cc}
0&1\\
1&0\\
\end{array}
\right]\in\pgl{2}{A}$ and for $u\in A^\times$ let $M_u:=M_{\infty,0,u}=\mathrm{diag}(u,1)$.

\begin{prop}\cite[Lemmas 4.4,4.7, Corolllaries 4.5, 4.8]{hut:f2} \label{prop:pres}
Let $A$ be a ring.
\begin{enumerate}
\item The $\pgl{2}{A}$-module $Z_2(\Sigma_\bullet)$ has the following presentation: It is generated by $\cconst{A}$ and $\suss{1}{u}$, $u\in A^\times$ subject to
the relations
\begin{enumerate}
\item $\cconst{A}\cdot\omega =\cconst{A}$,
\item $\cconst{A}\cdot\left(I+M_{0,1,\infty}+M_{1,\infty,0}\right)=0$,
\item $\suss{1}{u}\cdot\omega =\suss{1}{u^{-1}}$ for all $u\in A^\times$,
\item $\suss{1}{uv}=\suss{1}{u}+\suss{1}{v}\cdot M_u$ for all $u,v\in A^\times$.
\end{enumerate}
\item The $\sgr{A}$-module $\crpb{A}$ has the following presentation: It is generated $\cconst{A}$ and $\suss{1}{u}$, $u\in A^\times$ subject to the relations
\begin{enumerate}
\item $\an{-1}\cconst{A}=\cconst{A}$,
\item $3\cconst{A}=0$,
\item $\an{-1}\suss{1}{u}=\suss{1}{u^{-1}}$ for $u\in A^\times$,
\item $\suss{i}{uv}=\suss{i}{u}+\an{u}\suss{i}{v}$ for all $u,v\in A^\times$.
\end{enumerate}
\end{enumerate}
\end{prop}

\section{$3$-dimensional complexes}\label{sec:3dim}
Our main concern in this paper are local rings with residue field $\F{3}$. For such rings, the associated simplicial complex $Y(A)$ is $3$-dimensional. Recall that we are primarily interested in the $\pgl{2}{A}$-module $Z_2(L_\bullet)$ and its group of coinvariants.

We first investigate $Z_2(L_\bullet)$ for general $3$-dimensional simplicial complexes $Y$.  As we show below, understanding this group requires knowledge of the low-dimensional homology of the space $|Y|$ (Theorem \ref{thm:c3} and Corolllary \ref{cor:c3}). The main result of this section is Proposition \ref{prop:z3sigma}, which gives a complete family of generators for the group $Z_2(\Sigma_\bullet)\subset Z_2(L_\bullet)$, which in turn plays a key role in the results of the next section.

\subsection{Calculating $Z_2(L_\bullet)$ for $3$-dimensional complexes}

\begin{lem}\label{lem:c31} 
Let $Y$ be a simplicial complex with the property that every $2$-simplex is contained in a $3$-simplex. Then
\[
\bar{C}_3=d_4(\bar{C}_4)+\Z[\{ (x,x,x,y)\ | \{ x,y\}\in Y_1\}]+\Z[\{ (x,y,y,y)\ | \{ x,y\}\in Y_1\}]
+\Z[\{ (z,x,y,x)\  | \{ x,y,z\}\in Y_2\}].
\]
\end{lem}
\begin{proof} 
Let
\[
M:=d_4(\bar{C}_4)+\Z[\{ (x,x,x,y)\ | \{ x,y\}\in Y_1\}]+\Z[\{ (x,y,y,y)\ | \{ x,y\}\in Y_1\}]
+\Z[\{ (z,x,y,x)\  | \{ x,y,z\}\in Y_2\}].
\]
Given $\{ x,y,z\}\in Y_2$ we have $d_4(x,x,x,y,z)=(x,x,y,z)-(x,x,x,z)+(x,x,x,y)$ and thus $(x,x,y,z)\in M$. Similarly,
$d_4(x,y,z,z,z)=(y,z,z,z)-(x,z,z,z)+(x,y,z,z)$ and thus $(x,y,z,z)\in M$.  

Note that for all $x\in Y_0$, $d_4(x,x,x,x,x)=(x,x,x,x)$ and thus $(x,x,x,x)\in M$. Thus for all $\{ x,y\}\in Y_1$, $d_4(x,x,x,y,x)=(x,x,y,x)-(x,x,x,x)+(x,x,x,y)$ implies that $(x,x,y,x)\in M$ and likewise $d_4(x,y,x,x,x)=(y,x,x,x)-(x,x,x,x)+(x,y,x,x)$ implies that $(x,y,x,x)\in M$.

For all $\{ x,y,z\}\in Y_2$,  $d_4(x,y,y,z,y)=(y,y,z,y)-(x,y,y,y)+(x,y,y,z)$  and thus $(x,y,y,z)\in M$. 

Given $\{ x,y,z\}\in Y_2$, by hypothesis there exists $w\in Y_0$ with $\{ x,y,z,w\}\in Y_3$. Then $d_4(w,x,y,x,z)=(x,y,x,z)+(w,x,x,z)+(w,x,y,x)$ and hence $(x,y,x,z)\in M$.
Furthermore, $d_4(w,x,y,z,x)=(x,y,z,x)+(w,x,z,x)-(w,x,y,x)+(x,y,z,x)$ and hence $(x,y,z,x)\in M$. This proves the Lemma.
\end{proof}

\begin{lem}\label{lem:c32}  Let $Y$ be a simplicial complex satisfying 
condition $(C)_1$.

 Let $M_0\subset \bar{C}_3$ be the subgroup
\[
M_0:=d_4(\bar{C}_4)+\Z[\{ (x,x,x,y)\ | \{ x,y\}\in Y_1\}]+\Z[\{ (x,y,y,y)\ | \{ x,y\}\in Y_1\}].
\]
For any $\{ x,y\}\in Y_1$ and for any $z,w\in Y^{\{ x,y\}}_0$ we have
 $(z,x,y,x)\equiv (w,x,y,x)\mod{M_0}$.
\end{lem}

\begin{proof} Given $\{ x,y\}\in Y_1$ and $\{ z_1,z_2\}\in Y^{\{ x,y\}}_1$ we have $d(z_1,z_2,x,y,x)=(z_2,x,y,x)-(z_1,x,y,x)-(z_1,z_2,x,x)$. By the proof of Lemma \ref{lem:c31} (third sentence),
$(z_1,z_2,x,x)\in M_0$. Thus $(z_1,x,y,x)\equiv (z_2,x,y,x)\mod{M_0}$.

 Now let $z,w\in Y^{\{ x,y\}}_0$. Since $Y$ satisfies $(C)_1$, there exists an $\{ x,y\}$-path 
$z=z_0, z_1,\ldots,z_n=w$ in $Y$ and hence 
\[
(z,x,y,x)\equiv (z_1,x,y,x)\equiv \cdots \equiv (z_{n-1},x,y,x)\equiv (w,x,y,x)\mod{M_0}
\]
as claimed.
\end{proof}

\begin{prop}\label{prop:c3} Suppose that the simplicial complex $Y$ satisfies
\begin{enumerate}
\item every $1$-simplex is contained in a $2$-simplex, 
\item every $2$-simplex is contained in a $3$-simplex, and
\item condition $(C)_1$; i.e.,  for every $\alpha\in Y_1$, $Y^\alpha$ is path-connected.
\end{enumerate}
Then $H_3(\bar{C}_\bullet)=0$.
\end{prop}

\begin{proof} Let $h\in \bar{C}_3$. For each $\{ x,y\}\in Y_1$, fix $z_{x,y}\in Y^{ \{ x,y\}}_0$. By Lemmas \ref{lem:c31} and \ref{lem:c32}, we have 
\[
h=\sum_{(x,y)\in X_1}n_{x,y}(z_{x,y},x,y,x)+\sum_{(x,y)\in X_1}m_{x,y}(x,y,y,y)+\sum_{(x,y)\in X_1}\ell_{x,y}(x,x,x,y)+d_4(g)
\]
for some $g\in \bar{C}_4$. Thus, if $h\in Z_3(\bar{C}_\bullet)$ then
\begin{eqnarray*}
0=d_3(h)&=&\sum n_{x,y}(x,y,x)-\sum n_{x,y}(z_{x,y},x,x)\\
&&+\sum m_{x,y}(y,y,y)-\sum m_{x,y}(x,y,y)\\
&&+\sum\ell_{x,y}(x,x,y)-\sum{\ell_{x,y}}(x,x,x).
\end{eqnarray*}
By linear independence of terms, we must have $\sum n_{x,y}(x,y,x)=0$ and hence $n_{x,y}=0$ for all $(x,y)\in X_1$. Therefore the terms $\sum m_{x,y}(x,y,y)$ and 
$\sum\ell_{x,y}(x,x,y)$ must vanish and hence $m_{x,y}=0=\ell_{x,y}$ for all $(x,y)$. Thus $h=d_4(g)$, as required.
\end{proof}

\begin{thm}\label{thm:c3}  Let $Y$ be a simplicial complex such that $Y_4=\emptyset$ and suppose that $Y$ satisfies conditions (1)-(3) of Proposition \ref{prop:c3}. Then the natural map $H_2(L_\bullet)\to H_2(|Y|)$ is an isomorphism and there is a natural
exact sequence
\[
\xymatrix{
H_3(|Y|)\ar[r]& L_3/Z_3(\Sigma_\bullet)\ar^-{d_3}[r] &Z_2(L_\bullet)\ar[r]& H_2(|Y|)\ar[r]& 0.
}
\]
\end{thm}

\begin{proof} The homology long exact sequence associated to the short exact sequence of complexes $0\to L_\bullet\to C_\bullet\to \bar{C}_\bullet\to 0$ gives an exact sequence
$H_3(\bar{C}_\bullet)\to H_2(L_\bullet)\to H_2(|Y|)\to H_2(\bar{C}_\bullet)$. Since the first and last terms are $0$ by Lemma \ref{lem:c2} and Proposition \ref{prop:c3},
it follows that $H_2(L_\bullet)\cong H_2(|Y|)$. 

Now by definition, we have an exact sequence 
\[
0\to Z_3(L_\bullet)\to L_3\to Z_2(L_\bullet)\to H_2(L_\bullet)\to 0.
\]
Since $Z_3(\Sigma_\bullet)\subset Z_3(L_\bullet)$, we derive an exact sequence
\[
0\to Z_3(L_\bullet)/Z_3(\Sigma_\bullet)\to L_3/Z_3(\Sigma_\bullet)\to Z_2(L_\bullet)\to H_2(|Y|)\to 0.
\]

However, since $Y_4=\emptyset$, we have $L_4=0$ and hence $H_3(L_\bullet)=Z_3(L_\bullet)$. Thus we have an exact sequence
\[
\xymatrix{
H_4(\bar{C}_\bullet)\ar^-{\delta}[r]& Z_3(L_\bullet)\ar[r]&H_3(|Y|)\ar[r]& 0=H_3(\bar{C}_\bullet).
}
\]

But $\image{\delta}=d_4(D_4)\cap L_3$ is contained in $Z_3(\Sigma_\bullet)$ by Proposition \ref{prop:symm} and thus there is an induced surjective homomorphism
$\xymatrix{H_3(|Y|)\ar@{>>}[r]& Z_3(L_\bullet)/Z_3(\Sigma_\bullet)}$.
\end{proof}

\begin{cor}\label{cor:c3}  Let $Y$ be a simplicial complex such that $Y_4=\emptyset$ and suppose that $Y$ satisfies conditions (1)-(3) of Proposition \ref{prop:c3}. Suppose further that 
$H_2(|Y|)=0$. Then  there is a natural
exact sequence
\[
\xymatrix{
H_3(|Y|)\ar[r]& L_3/Z_3(\Sigma_\bullet)\ar^-{d_3}[r] &Z_2(L_\bullet)\ar[r]& 0.
}
\]
\end{cor}

Note that, by Proposition \ref{prop:localhom},  Corollary \ref{cor:c3} applies in particular to the simplicial complex $Y(A)$ when $A$ is a local ring with residue field $\F{3}$.
\subsection{Generators for $Z_3(\Sigma_\bullet)$}
Our objective is to give a presentation of $Z_2(L_\bullet)$, at least when $H_2(|Y|)=0$. Since we already have set of generators of  $L_3$ from its definition, in view of Theorem \ref{thm:c3}, this  motivates us to determine a set of generators for $Z_3(\Sigma_\bullet)$:

 Recall that $\Sigma_3$ is  generated by the elements $R_i^\alpha(a,b)\in L_3$ with $i\in \{ 1,2,3\}$, $\alpha:=\{ x,y\}\in Y_1$ and 
$\alpha\cup \{ a,b\}\in Y_3$. These elements satisfy
\begin{eqnarray*}
d_3(R_1^\alpha(a,b))&=&R_1^\alpha(b)-R_1^\alpha(a)=D_1^\alpha(a,b)\\
d_3(R_2^\alpha(a,b))&=&R_1^\alpha(b)-R_2^\alpha(a)=D_1^\alpha(a,b)+C^\alpha(a)\\
d_3(R_3^\alpha(a,b))&=&R_2^\alpha(b)-R_2^\alpha(a)=D_2^\alpha(a,b)=D_1^\alpha(a,b)+C^\alpha(a)-C^\alpha(b)\\
\end{eqnarray*}
Now let 
\begin{eqnarray*}
D_{12}^\alpha(a,b):=&R_2^\alpha(a,b)-R_1^\alpha(a,b)\\
D_{32}^\alpha(a,b):=&R_2^\alpha(a,b)-R_3^\alpha(a,b).\\
\end{eqnarray*}
Then we have
\begin{eqnarray*}
d_3(D_{12}^\alpha(a,b))&=& R_1^\alpha(a)-R_2^\alpha(a)=C^\alpha(a),\\
d_3(D_{32}^\alpha(a,b))&=& R_1^\alpha(b)-R_2^\alpha(b)=C^\alpha(b).\\
\end{eqnarray*}
Thus we have the following families of elements of $Z_3(\Sigma_\bullet)$:  
\begin{eqnarray*}
DD_1^\alpha(a,b,b'):=& D_{12}^\alpha(a,b)-D_{12}^\alpha(a,b'),\\
DD_2^\alpha(a,a',b):=& D_{32}^\alpha(a,b)-D_{12}^\alpha(b,a'),\\
\end{eqnarray*}
for $\alpha\in Y_1$, and $\{ a,b\}, \{ a,b'\},\{ a',b\}\in Y_1^\alpha$.

For $\alpha\in Y_1$ and if $\mathbf{p}=a_0,a_1,\ldots,a_n$ is an $\alpha$-path in $Y$, we will let
\[
R_i^\alpha(\mathbf{p}):=\sum_{i=1}^nR_i^\alpha(a_{i-1},a_i)\in \Sigma_3
\]
for $i=1,2$.

Observe the that $d_3(R_1^\alpha(\mathbf{p}))=\sum_{i=1}^n\left( R_1^\alpha(a_i)-R_1^\alpha(a_{i-1})\right)=R_1^\alpha(a_n)-R_1^\alpha(a_0)$.
Thus whenever $\mathbf{p}$ is an $\alpha$-loop $R_1^\alpha(\mathbf{p})\in Z_3(\Sigma_\bullet)$.

In this section we will prove:

\begin{prop}\label{prop:z3sigma}
Suppose that $Y$ satisfies condition $(C)_1$ and suppose also that every $2$-simplex is contained in a $3$-simplex. Then $Z_3(\Sigma_\bullet)$ is generated by the following three families:
\begin{enumerate}
\item $\{ {DD}_1^\alpha(a,b,b') \ |\ \alpha\in Y_1$, $\{ a,b\},\{ a,b'\}\in Y_1^\alpha\}$.
\item $\{ {DD}_2^\alpha(a,a',b) \ |\ \alpha\in Y_1$, $\{ a,b\},\{ a',b\}\in Y_1^\alpha\}$.
\item $\{ {R}_1^\alpha(\mathbf{p})\ |\ \alpha\in Y_1$, $\mathbf{p}\mbox{ is an $\alpha$-loop}\}$.
\end{enumerate}
\end{prop}

\hfill

Let $\widehat{\Sigma}_3$ be the free $\Z$-module with basis $\widehat{R}_i^\alpha(a,b)$,  with $i\in \{ 1,2,3\}$, $\alpha:=\{ x,y\}\in Y_1$ and 
$\alpha\cup \{ a,b\}\in Y_3$. 

Let $\widehat{Z}_2:=\widehat{Z}^1\oplus\widehat{Z}^c$ where $\widehat{Z}^1$ is the free $\Z$-module with basis the symbols $\widehat{D}_1^\alpha(a,b)$ for all $\alpha,a,b$ satisfying $\alpha\in Y_1$ and $a,b\in Y_0^\alpha$ and $\widehat{Z}^c$ is the free $\Z$-module with generators $\widehat{C}^\alpha(a)$, where  ${ a}\in Y_0^\alpha$. 

Let $ \widehat{d}_3:\widehat{\Sigma}_3\to \widehat{Z}_2$ be the homomorphism defined by 
\begin{eqnarray*}
\widehat{d}_3(\widehat{R}_1^\alpha(a,b))&=& \widehat{D}_1^\alpha(a,b),\\
\widehat{d}_3(\widehat{R}_2^\alpha(a,b))&=& \widehat{D}_1^\alpha(a,b)+\widehat{C}^\alpha(a),\\
\widehat{d}_3(\widehat{R}_3^\alpha(a,b))&=& \widehat{D}_1^\alpha(a,b)+\widehat{C}^\alpha(a)-\widehat{C}^\alpha(b).\\
\end{eqnarray*}
If we denote by $\pi_3:\widehat{\Sigma}_3\to\Sigma_3$ and $\pi_2=\pi_2^1\oplus\pi_2^c:\widehat{Z}_2\to Z_2(\Sigma_\bullet)=Z^1\oplus Z^c$ the obvious surjective homomorphisms then the diagram 
(with exact rows)
\[
\xymatrix{0\ar[r]&\ker{\widehat{d}_3}\ar[r]\ar^-{\pi_3}[d]&\widehat{\Sigma}_3\ar^-{\pi_3}[d]\ar^-{\widehat{d}_3}[r]\ar^-{\pi_2\circ\widehat{d}_3}[dr]&  \widehat{Z}_2=\widehat{Z}^1\oplus \widehat{Z}^c\ar^-{\pi_2=\pi_2^1\oplus\pi_2^c}[d]\\
0\ar[r]& Z_3(\Sigma_\bullet)\ar[r]&\Sigma_3\ar^-{d_3}[r]&Z_2(\Sigma_\bullet)=Z^1\oplus Z^c\\
}
\]
commutes.

Now let $\widehat{D}_{12}^\alpha(a,b):=\widehat{R}_2^\alpha(a,b)-\widehat{R}_1^\alpha(a,b)$ in $\widehat{\Sigma}_3$ and let 
 $\widehat{D}_{32}^\alpha(a,b):=\widehat{R}_2^\alpha(a,b)-\widehat{R}_3^\alpha(a,b)$.
Finally, in $\widehat{\Sigma}_3$, let 
\begin{eqnarray*}
\widehat{DD}_1^\alpha(a,b,b'):=& \widehat{D}_{12}^\alpha(a,b)-\widehat{D}_{12}^\alpha(a,b'),\\
\widehat{DD}_2^\alpha(a,a',b):=& \widehat{D}_{32}^\alpha(a,b)-\widehat{D}_{12}^\alpha(b,a'),\\
\end{eqnarray*}
\begin{lem}\label{lem:kerd3}
The $\Z$-module $\ker{\widehat{d}_3}$ is generated by the families $\{ \widehat{DD}_1^\alpha(a,b,b') \}$ and $\{ \widehat{DD}_2(a,a',b)\}$.
\end{lem}

\begin{proof} We begin by noting that there are obvious decompositions $\widehat{\Sigma}_3=\oplus_{\alpha\in Y_1}\widehat{\Sigma}_3^\alpha$, where $\widehat{\Sigma}_3^\alpha$ 
has basis $\{ \widehat{R}_i^\alpha(a,b)\}$ with $\alpha$ fixed and $\{ a,b\}$ satisfying $\alpha\cup \{ a,b\}\in Y^3$.  Likewise, $\widehat{Z}^1=\oplus_\alpha \widehat{Z}^{1,\alpha}$, $\widehat{Z}^c=\oplus_\alpha\widehat{Z}^{c,\alpha}$, where $\widehat{Z}^{1,\alpha}$ has basis $\{ \widehat{D}_1^\alpha(a,b)\}$ with $a, b$ satisfying $\{ a\}\cup\alpha,\{ b\}\cup\alpha\in Y_2$. Furthermore, $\widehat{d}_3=\oplus_\alpha \widehat{d}_3^\alpha$ where $\widehat{d}_3^\alpha:\widehat{\Sigma}_3\to \widehat{Z}_2^\alpha=\widehat{Z}^{1,\alpha}\oplus 
\widehat{Z}^{c,\alpha}$. Thus it is enough to prove that, for any $\alpha\in Y_1$,  $\ker{\widehat{d}_3^\alpha}$ is generated by the given families. 

Thus we fix $\alpha\in Y_1$ and prove this. (Since $\alpha$ is now fixed, for simplicity it is removed from the notation in the remainder of the proof.)

Now suppose that $w=\sum_{i=1}^3\left(\sum_{(a,b)}m^i_{a,b}\widehat{R}_i(a,b)\right)\in \ker{\widehat{d}_3}$. Here, the sum is over all pairs $(a,b)$ such that $\{ a,b\}\cup \alpha\in Y_3$, and $m^i_{a,b}\in \Z$ for all $i\in\{ 1,2,3\}$ and for all $(a,b)$.

Thus 
\begin{eqnarray*}
0=\widehat{d}_3(w)&=& \sum_{(a,b)}m^1_{a,b}\widehat{D}_1(a,b)+\sum_{(a,b)}m^2_{a,b}\left(\widehat{D}_1(a,b)+\widehat{C}(a)\right)\\
&&+\sum_{(a,b)}m^3_{a,b}\left( \widehat{D}_1(a,b)+\widehat(C)(a)-\widehat{C}(b)\right)\\
&=& \sum_{(a,b)}\left( m^1_{a,b}+m^2_{a,b}+m^3_{a,b}\right)\widehat{D}_1(a,b)\\
&&+ \sum_{a}\left(  \sum_b(m^2_{a,b}+m^3_{a,b})-\sum_bm^3_{b,a}\right)\widehat{C}(a).\\
\end{eqnarray*}
Thus it follows that both of these latter terms vanish in $\widehat{Z}_2$. The vanishing of the first gives the identity\\
 (i) $m^2_{a,b}=-m^1_{a,b}-m^3_{a,b}$ for all $(a,b)$.\\
Using this, the vanishing of the second then implies the identities\\
(ii) $\sum_bm^1_{a,b}=-\sum_bm^3_{b,a}$ for all $a$.

Now for each $a$ occurring, fix $t_a$ such $\{ a,t_a\}\cup \alpha\in Y_3$. (If no such $t_a$ exists, then there cannot be any symbols  $R_i(a,b)$ and $a$ cannot occur.)
Thus
\begin{eqnarray*}
-w&=&\sum_{(a,b)}m^1_{a,b}\left( \widehat{R}_2(a,b)-\widehat{R}_1(a,b)\right)+\sum_{(a,b)}m^3_{a,b}\left(\widehat{R}_2(a,b)-\widehat{R}_3(a,b)\right)\mbox{ (using (i)}\\
&=&\sum_{(a,b)}m^1_{a,b}\widehat{D}_{12}(a,b)+\sum_{(a,b)}m^3_{a,b}\widehat{D}_{32}(a,b)\\
&=& \sum_{(a,b)}m^1_{a,b}\left(\widehat{D}_{12}(a,b)-\widehat{D}_{12}(a,t_a)\right)+\sum_{(a,b)}m^1_{a,b}\widehat{D}_{12}(a,t_a)+\sum_{(a,b)}m^3_{a,b}\widehat{D}_{32}(a,b)\\
&=& \sum_{(a,b)}m^1_{a,b}\widehat{DD}_{1}(a,b,t_a)+\sum_{a}\left( \sum_b m^1_{a,b}\right)\widehat{D}_{12}(a,t_a)+\sum_{(a,b)}m^3_{a,b}\widehat{D}_{32}(a,b).\\
\end{eqnarray*}
However, by (ii), $ \sum_b m^1_{a,b}=-\sum_b m^3_{b,a}$ for all $a$. So
\begin{eqnarray*}
-w&=& \sum_{(a,b)}m^1_{a,b}\widehat{DD}_{1}(a,b,t_a)-\sum_{a}\left( \sum_b m^3_{b,a}\right)\widehat{D}_{12}(a,t_a)+\sum_{(b,a)}m^3_{b,a}\widehat{D}_{32}(b,a)\\
&=& \sum_{(a,b)}m^1_{a,b}\widehat{DD}_{1}(a,b,t_a)+\sum_{(b,a)}m^3_{b,a}\left( \widehat{D}_{32}(b,a)-\widehat{D}_{12}(a,t_a)\right)\\
&=& \sum_{(a,b)}m^1_{a,b}\widehat{DD}_{1}(a,b,t_a)+\sum_{(b,a)}m^3_{b,a}\widehat{DD}_2(b, t_a,a)
\end{eqnarray*}
as required.
\end{proof}

Now let $\alpha=\{ x,y\}\in Y_1$. Recall that $\widehat{\Sigma}_3^\alpha$ is free with basis $\{ \widehat{R}_i^\alpha(a,b)\  |\ 1\leq i\leq 3, \{ a,b\}\in Y^\alpha_1\}$ while $\widehat{Z}^{1,\alpha}$ is free with basis $\{ \widehat{D}_1^\alpha(a,b)\ |\ a,b\in Y_0^\alpha\}$. 

We let $\widehat{Z}^{1,\alpha}_w$ denote the free subgroup of $\widehat{Z}^{1,\alpha}$ with basis $\{ \widehat{D}_1^\alpha(a,b)\ |\ \{ a, b\}\in Y_1^\alpha\}$, and we let $\widehat{Z}^1_w:=\oplus_{\alpha\in Y_1}\widehat{Z}^{1,\alpha}_w\subset \widehat{Z}^1$.

\begin{lem}\label{lem:image} Suppose that $Y$ satisfies the condition that for any $\alpha\in Y_1$ and for any $a\in Y^\alpha_0$, there exists $b\in Y^\alpha_0$ with $\{ a,b\}\in Y^\alpha_1$. Then $\image{\widehat{d}_3}$ is $\widehat{Z}^1_w\oplus \widehat{Z}^c\subset \widehat{Z}^1\oplus\widehat{Z}^c=\widehat{Z}_2$.
\end{lem}

\begin{proof} This follows immediately from the definition of $\widehat{\Sigma}_3$ and $\widehat{d}_3$, together with the fact that $\widehat{Z}^c\subset \image{\widehat{d}_3}$ 
since $\widehat{C}^\alpha(a)=\widehat{d}_3(\widehat{D}_{12}^\alpha(a,b))$  when $\{ a,b\}\in Y^\alpha_1$.
\end{proof}

By Proposition \ref{cor:z2pres}, 
$\ker{\pi_2^c}$ is generated by the family $\{ \widehat{C}^{x,y}(z)+\widehat{C}^{x,z}(y)+\widehat{C}^{y,z}(x)\}\subset \widehat{Z}^c$ and $\ker{\pi_2^1}$ 
is generated by the family
$\{ \widehat{D}_1^\alpha(a,b)+\widehat{D}_1^\alpha(b,c)-\widehat{D}_1^\alpha(a,c)\ | \alpha\in Y_1,\ a,b,c\in Y^\alpha_0\}\subset\widehat{Z}^1$. (Recall that here we allow $a=c$ etc and we define $\widehat{D}_1^\alpha(a,a):=0$ in this case.)

If $\alpha\in Y_1$ and if $\mathbf{p}=a_0,a_1,\ldots,a_n$ is an $\alpha$-path in $Y$, we will let
\[
\widehat{D}_1^\alpha(\mathbf{p}):=\sum_{i=1}^n\widehat{D}_1^\alpha(a_{i-1},a_i)\in \widehat{Z}^{1,\alpha}_w\subset \widehat{Z}^1.
\]
\begin{lem} \label{lem:ker1} Suppose that $Y$ satisfies condition $(C)_1$. Then $\widehat{Z}^1_w\cap \ker{\pi_2^1}$ is generated by $\{ \widehat{D}_1^\alpha(\mathbf{p})\ |\ \alpha\in Y_1, \mathbf{p}
\mbox{ is an $\alpha$-loop in } Y\}$.
\end{lem}

\begin{proof} For the proof, we may fix $\alpha\in Y_1$ and remove references to $\alpha$ from our notation. We let $K_w$ be the subgroup of $\widehat{Z}^1_w$ generated by 
the set   $\{ \widehat{D}_1(\mathbf{p})\ |\ \mathbf{p}
\mbox{ is an $\alpha$-loop in } Y\}$.  

We begin by observing that if $\mathbf{p}=a_0,a_1,\ldots,a_n$ is an $\alpha$-loop then $\widehat{D}_1(\mathbf{p})\in \ker{\pi_2^1}$; i.e., $\widehat{D}_1(\mathbf{p})\equiv 0\mod{\ker{\pi_2^1}}$: For all $i$, we have $\widehat{D}_1(a_{i-1},a_i)+\widehat{D}_1(a_i,a_{i+1})\equiv \widehat{D}_1(a_{i-1},a_{i+1})\mod{\ker{\pi_2^1}}$ and thus, by induction, $\widehat{D}_1(\mathbf{p})\equiv\widehat{D}_1(a_0,a_n)\equiv 0\mod{\ker{\pi_2^1}}$ since $a_0=a_n$. Thus, we have $K_w\subset\ker{\pi_2^1}$.

If $\mathbf{p}=a_0,\ldots,a_n$ and $\mathbf{q}=b_0,\ldots,b_m$ are $\alpha$-paths satsifying $a_n=b_0$, we let $\mathbf{p}\ast\mathbf{q}$ be the $\alpha$-path 
$a_0,\ldots,a_n-b_0, b_1,\ldots,b_m$ and note that $\widehat{D}_1(\mathbf{p}\ast\mathbf{q})=\widehat{D}_1(\mathbf{p})+\widehat{D}_1(\mathbf{q})$. Futhermore, 
let $\bar{\mathbf{p}}$ denote the $\alpha$-path $a_n,\ldots,a_0$ and observe that $\widehat{D}_1(\mathbf{p})+\widehat{D}_1(\bar{\mathbf{p}})=\widehat{D}_1(\mathbf{p}\ast\bar{\mathbf{p}})\in K_w$, so that $\widehat{D}_1(\bar{\mathbf{p}})\equiv -\widehat{D}_1(\mathbf{p})\mod{K_w}$.

Given $a,b\in Y^\alpha_0$ let $\mathbf{p}=a_0,\ldots,a_n$ and $\mathbf{q}=b_0,\ldots,b_m$ be any $\alpha$-paths from $a$ to $b$. Then $\mathbf{p}\ast\bar{\mathbf{q}}$ is an $\alpha$-loop at $a$. Hence
\[
\widehat{D}_1(\mathbf{p}\ast\bar{\mathbf{q}})=\widehat{D}_1(\mathbf{p})+\widehat{D}_1(\bar{\mathbf{q}})=\widehat{D}_1(\mathbf{p})-\widehat{D}_1(\mathbf{q})\in K_w
\]
and hence $\widehat{D}_1(\mathbf{p})\equiv\widehat{D}_1(\mathbf{q})\mod{K_w}$. Thus the following homomorphism is well-defined:\\
 $f:\widehat{Z}^1\to \widehat{Z}^1_w/K_w$  sends $\widehat{D}_1(a,b)$ to 
$\widehat{D}_1(\mathbf{p})$, where $\mathbf{p}$ is any $\alpha$-path from $a$ to $b$. Furthermore, if $\mathbf{p}$ is an $\alpha$-path from $a$ to $b$, and if $\mathbf{q}$ is an $\alpha$-path from $b$ to $c$, then $\mathbf{p}\ast\mathbf{q}$ is an $\alpha$-path from $a$ to $c$ and hence
\[
f(\widehat{D}_1(a,c))=\widehat{D}_1(\mathbf{p}\ast\mathbf{q})=\widehat{D}_1(\mathbf{p})+\widehat{D}_1(\mathbf{q})=f(\widehat{D}_1(a,b)+\widehat{D}_1(b,c)).
\]
Hence $\ker{\pi_2^1}\subset \ker{f}$. Thus $f$ induces a homomorphism $\widehat{Z}^1/\ker{\pi_2^1}\to \widehat{Z}^1_w/K_w$ which is inverse to the natural map $\widehat{Z}^1_w/K_w\to \widehat{Z}^1/\ker{\pi_2^1}$ induced by the inclusion $\iota:\widehat{Z}^1_w\to \widehat{Z}^1$. Thus $\iota$ induces isomorphisms of groups
\[
\xymatrix{
{\widehat{Z}^1_w}/{K_w}\ar[r]^-{\cong}&{\widehat{Z}^1}/{\ker{\pi_2^1}}\ar[r]^-{\cong}&  Z^1.\\
}
\]
In particular, the left map is injective and so the Lemma is proved.
\end{proof}

Lemmas  \ref{lem:image}, \ref{lem:ker1} and Proposition \ref{cor:z2pres} now immediately imply
\begin{cor} \label{cor:imker}
 Suppose that $Y$ satisfies condition $(C)_1$. Then $\image{\widehat{d}_3}\cap\ker{\pi_2}$ is generated by the families of elements\\
  $\{ \widehat{D}_1^\alpha(\mathbf{p})\ |\ \alpha\in Y_1, \mathbf{p}
\mbox{ is an $\alpha$-loop in } Y\}$ and\\
 $\{ \widehat{C}^{\{ x,y\}}(z)+ \widehat{C}^{\{ x,z\}}(y)+\widehat{C}^{\{ y,z\}}(x)\ |\ \{x,y,z\}\in Y_2\}$.
\end{cor}

Recall again that when $\{ a,b\}\in Y_1^\alpha$, $\widehat{d}_3(\widehat{R}_1^\alpha(a,b))=\widehat{D}_1^\alpha(a,b)$. Hence if $\mathbf{p}$ is an $\alpha$-path in $Y$ then
$\widehat{d}_3(\widehat{R}_1^\alpha(\mathbf{p}))=\widehat{D}_1^\alpha(\mathbf{p})$.  Furthermore if $a\in Y_0^\alpha$ then $\widehat{d}_3(\widehat{D}_{12}^\alpha(a,b))=\widehat{C}^\alpha(a)$ for any $b$ such that $\{ a,b\}\in Y_1^\alpha$. Thus if $\{ x,y,z,b\}\in Y_3$ then
\[
\widehat{d}_3\left(\widehat{D}_{12}^{\{ x,y\}}(z,b))+\widehat{D}_{12}^{\{ x,z\}}(y,b))+\widehat{D}_{12}^{\{ y,z\}}(x,b))\right)=
 \widehat{C}^{\{ x,y\}}(z)+ \widehat{C}^{\{ x,z\}}(y)+\widehat{C}^{\{ y,z\}}(x).
\]

Thus, combining Lemma \ref{lem:kerd3} and Corollary \ref{cor:imker} we deduce:
\begin{lem}\label{lem:pi2d3}
Suppose that $Y$ satisfies condition $(C)_1$ and suppose also that every $2$-simplex is contained in a $3$-simplex. Then $\ker{\pi_2\circ \widehat{d}_3}=\ker{d_3\circ \pi_3}$ is generated by the following four families of elements:
\begin{enumerate}
\item $\{ \widehat{DD}_1^\alpha(a,b,b') \ |\ \alpha\in Y_1,\{ a,b\},\{ a,b'\}\in Y_1^\alpha\}$.
\item $\{ \widehat{DD}_2^\alpha(a,a',b) \ |\ \alpha\in Y_1, \{ a,b\},\{ a',b\}\in Y_1^\alpha\}$.
\item $\{ \widehat{R}_1^\alpha(\mathbf{p})\ |\ \alpha\in Y_1, \mathbf{p} \mbox{ is an $\alpha$-loop}\}$.
\item $\{ (\widehat{D}_{12}^{\{ x,y\}}(z,b))+\widehat{D}_{12}^{\{ x,z\}}(y,b))+\widehat{D}_{12}^{\{ y,z\}}(x,b))\}\ |\ \{ x,y,z,b\}\in Y_3\}$.
\end{enumerate}
\end{lem}

Noting that, since $\pi_3$ is surjective, we have  $\pi_3\left(\ker{d_3\circ \pi_3}\right)=\ker{d_3}=Z_3(\Sigma_\bullet)$, we conclude the proof of Proposition \ref{prop:z3sigma} by applying $\pi_3$:

\begin{proof}   For any $\{ x,y,z,b\}\in Y_3$, we have 
\[
\pi_3(\widehat{D}_{12}^{\{ x,y\}}(z,b))+\widehat{D}_{12}^{\{ x,z\}}(y,b))+\widehat{D}_{12}^{\{ y,z\}}(x,b)))=
{D}_{12}^{\{ x,y\}}(z,b))+{D}_{12}^{\{ x,z\}}(y,b))+{D}_{12}^{\{ y,z\}}(x,b))=0
\]
in $\Sigma_3$, so that the family (4) maps to $0$ in $Z_3(\Sigma_\bullet)$.
\end{proof}


\section{The module $\hrpb{A}$} \label{sec:hrpba}


Let $A$ be a ring satisfying (i) $L_4(A)=0$,  (ii) $\wn{A}$ generates $A^\times$ and (iii) $H_2(|Y(A)|)=0$  (eg. if $A$ is local with residue field $\F{3}$). By Corollary \ref{cor:c3} there is an exact sequence 
\[
H_3(|Y(A)|)\to\frac{L_3(A)}{Z_3(\Sigma_\bullet)}\to Z_2(L_\bullet)\to 0.
\]
Thus taking $\spl{2}{A}$-coinvariants we derive an exact sequence 
\[
H_3(|Y(A)|)\to\hrpb{A}\to\rpb{A}\to 0
\]
where $\hrpb{A}:=\left(L_3(A)/Z_3(\Sigma_\bullet)\right)_{\spl{2}{A}}$ (an $\sgr{A}$-module). 
The objective of this section is to give a presentation of the module $\hrpb{A}$ (Lemma \ref{lem:hrpbapres}). Lemma \ref{lem:hrpbapres} shows that $\hrpb{A}$ is generated by 
the elements $\Suss{1}{u}$  (which map to $\suss{i}{u}$ in $\rpb{A}$), $u\in \wn{A}$ and the constant $\cconst{A}$, subject to known generic families of relations or identities. 
\subsection{Basic generators}
Let $A$ be a ring. 
For $u\in \wn{A}$, we let $\gpb{u}:=(\infty,0,1,u)\in L_3(A)$. Thus, for $\{ a,b,c,d\}\in Y_3(A)$, we have 
\[
(a,b,c,d)=(\infty,0,1,d\cdot M_{a,b,c}^{-1})\cdot M_{a,b,c}=(\infty,0,1,u(a,b,c,d))\cdot M_{a,b,c}=\gpb{u(a,b,c,d)}M_{a,b,c}
\]
where 
\[
u(a,b,c,d):=d\cdot  M_{a,b,c}^{-1}=\frac{d(a,c)d(b,d)}{d(a,d)d(b,c)}\in \wn{A}.
\]

Since the stabilizer in $\pgl{2}{A}$ of $(\infty,0,1)$ is trivial, it follows that $L_3(A)$ is a free $\Z[\pgl{2}{A}]$-module with basis $\{ \gpb{u}\ |\ u\in \wn{A}\}$.

\subsection{The action of the symmetric group of degree $3$}
Let $S_3$ denote the group of permutations of $\{ 1,2,3\}$ acting on the right. Thus $S_3$ acts \emph{on the left} on the set $X_2(A)$, and hence on $L_2(A)$, by permuting the components:
$\sigma\cdot (a_1,a_2,a_3):= (a_{1\cdot\sigma},a_{2\cdot\sigma},a_{3\cdot\sigma})$. This action commutes with action of $\pgl{2}{A}$, which acts componentwise.
We let $S_3$ act on the left on $X_3(A)$, and hence on $L_3(A)$, via
\[
\sigma\ast(a_1,a_2,a_3,a_4):=( \sigma\cdot (a_1,a_2,a_3),a_4)=(a_{1\cdot\sigma},a_{2\cdot\sigma},a_{3\cdot\sigma},a_4).
\]
This action again commutes with the action of $\pgl{2}{A}$: For all $\ell\in L_3(A)$, $\gamma\in\Z[S_3]$, $M\in \pgl{2}{A}$ we have
\begin{eqnarray}\label{equivar}
\gamma\ast(\ell\cdot M)=(\gamma\ast \ell)\cdot M.
\end{eqnarray}

Furthermore, we fix the bijection between $\{ \infty,0,1\}$ and $\{ 1,2,3\}$ sending $\infty$ to $1$, $0$ to $2$ and $1$ to $3$, so that $S_3$ acts (on the right) as permutations of
$\{\infty,0,1\}$ via this bijection. For $\sigma\in S_3$ we now set
\[
M_\sigma:=M_{\infty\cdot\sigma,0\cdot\sigma,1\cdot\sigma}\in \pgl{2}{A}.
\]
This gives a group embedding $S_3\to \pgl{2}{A}$, $\sigma\mapsto M_\sigma$. We set $\tau:=(12)$,$\rho:=(123)\in S_3$. Thus
\[
M_\tau=\omega=
\left[
\begin{array}{cc}
0&1\\
1&0\\
\end{array}
\right],\quad 
M_{\rho}=M_{0,1,\infty}=
\left[
\begin{array}{cc}
0&-1\\
1&1\\
\end{array}
\right].
\]
We observe also that $\mathrm{det}(M_\sigma)=\mathrm{sgn}(\sigma)$ for all $\sigma\in S_3$.

With this definition, we have $\sigma\cdot (\infty,0,1)=(\infty\cdot\sigma,0\cdot\sigma,1\cdot\sigma)=(\infty,0,1)\cdot M_\sigma$ for all $\sigma\in S_3$.

Note then that for $u\in \wn{A}$, $\sigma\in S_3$, we have 
\begin{eqnarray*}
\sigma\ast\gpb{u}&=& \sigma\ast (\infty,0,1,u)
=(\sigma\cdot (\infty,0,1), u)\\
&=& ( (\infty,0,1)\cdot M_\sigma,u)=\left( \infty,0,1,u\cdot M_\sigma^{-1}\right)\cdot M_\sigma\\
&=&\gpb{u\cdot M_\sigma^{-1}}\cdot M_\sigma\\
&=& \gpb{\sigma\cdot u}\cdot M_\sigma
\end{eqnarray*}
where $\sigma\cdot u:= u\cdot M_\sigma^{-1}=u\cdot M_{\sigma^{-1}}$ for all $u\in \wn{A}$, $\sigma\in S_3$. For the purposes of reference below, we make this action of 
$S_3$ on $\wn{A}$ explicit: For all $u\in \wn{A}$, 
\begin{eqnarray*}
1\cdot u=u, & \rho\cdot u=\frac{u-1}{u},& \rho^2\cdot u =\frac{1}{1-u},\\
\tau\cdot u= \frac{1}{u},&(\tau\rho)\cdot u=\frac{u}{u-1},&(\tau\rho^2)\cdot u = 1-u.
\end{eqnarray*}

\subsection{Special elements and their properties}
We define the following elements in $\Sigma_3(A)\subset L_3(A)$:

For $u\in \wn{A}$, let 
\[
\Suss{1}{u}:=R_1^{\{\infty,0\}}(1,u)=(\infty,0,1,u)+(0,\infty,1,u)=(1+\tau)\ast\gpb{u}=\gpb{u}+\gpb{u^{-1}}\omega,
\]
and let
\begin{eqnarray*}
\Suss{2}{u}:&=&R_3^{\{\infty,0\}}(1,u)=(1,u,\infty,0)+(1,u,0,\infty)\\
&=&\gpb{u}\cdot
\left[
\begin{array}{cc}
-1&-1\\
u&1\\
\end{array}
\right]
+\gpb{u^{-1}}\cdot
\left[
\begin{array}{cc}
-1&-1\\
1&u^{-1}\\
\end{array}
\right].\\
\end{eqnarray*}
Note that $d_3(\Suss{i}{u})=\suss{i}{u}$ for $i\in \{ 1,2\}$ and for any $u\in \wn{A}$.

For  $u\in\wn{A}$, we define also the following element of $\Sigma_3(A)$:
\begin{eqnarray*}
\cconst{A}(u):&=&D^{\{ \infty,0\}}_{12}(1,u)=R_2^{\{ \infty,0\}}(1,u)-R_1^{\{ \infty,0\}}(1,u)\\
&=&(1,\infty,0,u)+(1,0,\infty,u)-\Suss{1}{u}\\
&=& 
\gpb{\frac{1}{1-u}}\cdot
M_{\rho^2}
+\gpb{\frac{1}{1-u^{-1}}}\cdot M_{\tau\rho}
-\Suss{1}{u}\\
&=&\Omega\ast\gpb{u}
\end{eqnarray*}
where $\Omega:=\rho^2+\tau\rho-1-\tau=(\rho^2-1)(1+\tau)\in \Z[S_3]$.

Note that $d_3(\cconst{A}(u))=C^{ \{ \infty,0\}}(1)=\cconst{A}$ for all $u\in \wn{A}$.

\begin{lem}\label{lem:cconstids} We have the following identities in $\Sigma_3(A)\subset L_3(A)$.
\begin{enumerate}
\item For all $u\in \wn{A}$ and for $i\in \{ 1,2\}$, we have $\Suss{i}{u}\cdot\omega=\Suss{i}{u^{-1}}$.
\item For all $u\in \wn{A}$, $\cconst{A}(u)\cdot \omega=\cconst{A}(u^{-1})$.
\item For all $u\in\wn{A}$,
\[
\cconst{A}(u)+\cconst{A}(1-u^{-1})\cdot M_{0,1,\infty}+\cconst{A}((1-u)^{-1})\cdot M_{1,\infty,0}=0.
\]
\end{enumerate}
\end{lem}

\begin{proof}
\begin{enumerate}
\item Note that $\infty\cdot\omega=0$, $0\cdot\omega=\infty$ and $u\cdot\omega=u^{-1}$ for all $u\in A^\times$.
Thus $\Suss{i}{u}\cdot\omega=R_j^{\{\infty,0\}}(1,u)\cdot\omega=R_j^{\{\infty,0\}}(1,u^{-1})=\Suss{i}{u^{-1}}$, as required.
\item Similarly, $\cconst{A}(u)\cdot\omega=D_{12}^{\{\infty,0\}}(1,u)\cdot\omega=D_{12}^{\{\infty,0\}}(1,u^{-1})=\cconst{A}(u^{-1})$.
\item  Note that for all $u\in \wn{A}$ and $\sigma\in S_3$, we have
\[
\cconst{A}(\sigma\cdot u)\cdot M_\sigma=(\Omega\ast\gpb{\sigma\cdot u})\cdot M_\sigma=\Omega\ast\left(\gpb{\sigma\cdot u}\cdot M_\sigma\right)
=\Omega\ast(\sigma\ast\gpb{u})=(\Omega\cdot\sigma)\ast\gpb{u}.
\]
Thus the given expression is $\left(\Omega\cdot (1+\rho+\rho^2)\right)\ast\gpb{u}$, which is $0$ since\\
 $\Omega\cdot (1+\rho+\rho^2)=(1+\rho+\rho^2)\cdot\Omega =0$ in $\Z[S_3]$.
\end{enumerate}
\end{proof}

\subsection{Presentation of the $\pgl{2}{A}$-module $L_3(A)/Z_3(\Sigma_\bullet)$}
Recall that for $u\in A^\times$, we let $M_u:=M_{\infty,0,u}=\mathrm{diag}(u,1)\in \pgl{2}{A}$.  The set $\{ M_u\ | u\in A^\times\}$ is clearly  a subgroup of $\pgl{2}{A}$ 
isomorphic to $A^\times$. The action of this group thus gives $L_n(A)$ the structure of a $\Z[A^\times]$-module via $\an{u}\cdot w:=w\cdot M_u$.

\begin{thm}\label{thm:l3pres} Let $A$ be a ring which satisfies the condition that $\wn{A}$ generates $A^\times$. Then $L_3(A)/Z_3(\Sigma_\bullet)$ has the following presentation as a $\pgl{2}{A}$-module:
It is generated by the elements $\gpb{u}$, $u\in\wn{A}$, subject to the relations
\begin{enumerate}
\item $\cconst{A}(u)=\cconst{A}(v)$ for all $u,v\in\wn{A}$,
\item $\an{u}\cconst{A}(u^{-1})-\cconst{A}(u)=\Suss{1}{u}-\Suss{2}{u}$ for all $u\in \wn{A}$, and
\item Whenever $u_1,\ldots,u_n\in \wn{A}$ satisfy $\prod_{i=1}^nu_i=1$ in $A^\times$ then
\[
\Sigma_{i=1}^n\an{u_1\cdots u_{i-1}}\Suss{1}{u_i}=\Suss{1}{u_1}+\an{u_1}\Suss{1}{u_2}+\cdots +\an{u_1\cdots u_{n-1}}\Suss{1}{u_n}=0.
\]
\end{enumerate}
\end{thm}

\begin{proof} Since $L_3(A)$ is a free $\pgl{2}{A}$-module with generators $\gpb{u}$, $u\in \wn{A}$, a set of generators of $Z_3(\Sigma_\bullet)$ give a complete set of relations 
among the elements $\gpb{u}$. By Propostion \ref{prop:z3sigma}, a set of generators of $Z_3(\Sigma_\bullet)$, and hence a set of relations of $L_3(A)/Z_3(\Sigma_\bullet)$,
 is the union of the following three families:
\begin{enumerate}
\item ${DD}_1^\alpha(a,b,b') =D_{12}^\alpha(a,b)-D_{12}^\alpha(a,b')$, $\alpha\in Y_1$, $\{ a,b\},\{ a,b'\}\in Y_1^\alpha$.
\item ${DD}_2^\alpha(a,a',b) =D_{32}^\alpha(a,b)-D_{12}^\alpha(b,a')$, $\alpha\in Y_1$, $\{ a,b\},\{ a',b\}\in Y_1^\alpha$.
\item $\{ {R}_1^\alpha(\mathbf{p})\}$, $\alpha\in Y_1$, $\mathbf{p}$ an $\alpha$-loop.
\end{enumerate}
However, note that ${DD}_2^\alpha(a,a',b) ={DD}_2^\alpha(a,a,b) +{DD}_1^\alpha(b,a,a')$, so that we can replace the second family with the family
\begin{enumerate}
\item[(2')] ${DD}_2^\alpha(a,a,b)=D_{32}^\alpha(a,b)-D_{12}^\alpha(b,a)$, $\alpha\in Y_1$, $\{ a,b\}\in Y_1^\alpha$.
\end{enumerate}

Let $\alpha=\{ x,y\}\in Y_1(A)$, and let $\{ a,b\}, \{ a,b'\}\in Y_1^\alpha$. Then 
\[
{DD}_1^\alpha(a,b,b')={DD}_1^{\{ \infty,0\}}(1,u,v)\cdot M_{x,y,a}
\]
where $u=u(x,y,a,b),v=u(x,y,a,b')\in \wn{A}$. Thus
\begin{eqnarray*}
{DD}_1^\alpha(a,b,b')&=&\left( D_{12}^{\{ \infty,0\}}(1,u)-D_{12}^{\{ \infty,0\}}(1,v)   \right)\cdot M_{x,y,a}\\
&=& \left(  \cconst{A}(u)-\cconst{A}(v)  \right)\cdot M_{x,y,a}.\\
\end{eqnarray*}
This accounts for the relation family (1) in the statement of the theorem.

Similarly,
\[
{DD}_2^\alpha(a,a,b)={DD}_2^{\{ \infty,0\}}(1,1,u)\cdot M_{x,y,a}
\]
where $u=u(x,y,a,b)$. Thus
\begin{eqnarray*}
{DD}_2^\alpha(a,a,b)&=&\left(  D_{32}^{\{ \infty,0\}}(1,u)-D_{12}^{\{ \infty,0\}}(u,1)\right)\cdot M_{x,y,a}\\
\end{eqnarray*}
However,
\begin{eqnarray*}
D_{32}^{\{ \infty,0\}}(1,u)&=& R_2^{\{ \infty,0\}}(1,u)-R_3^{\{ \infty,0\}}(1,u)\\
&=&\left(  R_2^{\{ \infty,0\}}(1,u)-R_1^{\{ \infty,0\}}(1,u)\right)+\left(  R_1^{\{ \infty,0\}}(1,u)-R_3^{\{ \infty,0\}}(1,u) \right)\\
&=&\cconst{A}(u)+\Suss{1}{u}-\Suss{2}{u}
\end{eqnarray*}
while
\[
D_{12}^{\{ \infty,0\}}(u,1)=D_{12}^{\{ \infty,0\}}(1,u^{-1})\cdot M_u=\cconst{A}(u^{-1})\cdot M_u.
\]
Hence
\[
{DD}_2^\alpha(a,a,b)=\left( \cconst{A}(u)+\Suss{1}{u}-\Suss{2}{u}-\cconst{A}(u^{-1})\cdot M_u \right)\cdot M_{x,y,a}
\]
and this gives us the family of relations (2) in the statement of the theorem.

Finally, let $\alpha=\{ x,y\}\in Y_1(A)$, and suppose that $\mathbf{p} =(a=a_0,a_1,\ldots, a-{n-1},a_n=a)$ is an $\alpha$-loop in $Y(A)$. Then
\begin{eqnarray*}
R_1^\alpha(\mathbf{p})&=&R_1^{ \{ x,y\}}(a_0,a_1)+R_1^{ \{ x,y\}}(a_1,a_2)+\cdots+R_1^{ \{ x,y\}}(a_{n-1},a_n)\\
&=&\left( R_1^{\infty,0}(1,w_1)+ R_1^{\infty,0}(w_1,w_2)+\cdots + R_1^{\infty,0}(w_{n-1},1)\right)\cdot M_{x,y,a}\\
\end{eqnarray*}
where $w_i=u(x,y,a,a_i)\in\wn{A}$.

 Since $1,w_1,\ldots, w_{n-1},w_n=1$ is an $\{ \infty,0\}$-loop in $Y(A)$,  by the proof of Lemma \ref{lem:ya} (3) there exist units $u_1,\ldots, u_n\in A^\times$ with $w_i=u_1\cdots u_i$ and $u_1\cdots  u_n=1$. Note also that
\[
 R_1^{\infty,0}(u,v)= R_1^{\infty,0}(1,v/u)M_u
\]
for any $u,v\in A^\times$. Thus 
\begin{eqnarray*}
R_1^{\infty,0}(1,w_1)+ R_1^{\infty,0}(w_1,w_2)+\cdots + R_1^{\infty,0}(w_{n-1},1)\\
=R_1^{\infty,0}(1,u_1)+ R_1^{\infty,0}(1,u_2)M_{w_1}+\cdots + R_1^{\infty,0}(1,u_n)M_{w_{n-1}}\\
=\Suss{1}{u_1}+\Suss{1}{u_2}\cdot M_{u_1}+\cdots +\Suss{1}{u_n}\cdot M_{u_1\cdots u_{n-1}}.
\end{eqnarray*}

Conversely, if $u_1,\ldots,u_n\in \wn{A}$ satisfy $\prod_{i=1}^nu_i=1$, then $\mathbf{p}=(1,w_1,\ldots, w_{n-1},w_n=1)$, where $w_i=u_1\cdots u_i$ is an $\{ \infty,0\}$-loop 
at $1$ in $Y(A)$ and 
\[
R_1^{\{ \infty,0\}}(\mathbf{p})=\Suss{1}{u_1}+\Suss{1}{u_2}\cdot M_{u_1}+\cdots +\Suss{1}{u_n}\cdot M_{u_1\cdots u_{n-1}}.
\]
Therefore generator family (3) of Proposition \ref{prop:z3sigma} gives the family of relations (3) of the statement of the theorem.
\end{proof}

\begin{cor}\label{cor:suss1}  Let $A$ be a ring which satisfies the condition that $\wn{A}$ generates $A^\times$.  
\begin{enumerate}
\item $\Suss{1}{u^{-1}}=-\an{u^{-1}}\Suss{1}{u}$ for all $u\in \wn{A}$.
\item 
Suppose that $u_1,\ldots,u_n$, $w_1,\ldots,w_m\in \wn{A}$ satisfy 
$u_1\cdots u_n=u=w_1\cdots w_m$. Then 
\[
\Sigma_{i=1}^n\an{u_1\cdots u_{i-1}}\Suss{1}{u_i}=\Sigma_{j=1}^m\an{w_1\cdots w_{j-1}}\Suss{1}{w_j}\mbox{ in }L_3(A)/Z_3(\Sigma_\bullet).
\]
\item By (2), for any $u\in A^\times$ we can define
\[
\Suss{1}{u}:=\Sigma_{i=1}^n\an{u_1\cdots u_{i-1}}\Suss{1}{u_i}\mbox{ in } L_3(A)/Z_3(\Sigma_\bullet)
\]
whenever $u_1,\ldots u_n$ satisfy $\prod_1^n u_i=u$. Then
\[
\Suss{1}{uv}=\Suss{1}{u}+\an{u}\Suss{1}{v}
\]
for all $u,v\in A^\times$.
\end{enumerate}
\end{cor}

\begin{proof}  
\begin{enumerate}
\item Taking $u_1=u^{-1},u_2=u$ in Theorem \ref{thm:l3pres} (3), we get 
\[
0=\Suss{1}{u^{-1}}+\an{u^{-1}}\Suss{1}{u}.
\]
\item Suppose that $\prod_1^n u_i=\prod_1^m w_j$ with $u_i,w_j\in \wn{A}$. Then $v_1,\ldots,v_{n+m}:=u_1,\ldots,u_n, w_m^{-1},\ldots,w_1^{-1}$ satisfy 
$\prod_1^{n+m}v_k=1$ and $v_k\in \wn{A}$ for all $k$. Thus by Theorem \ref{thm:l3pres} (3),
\begin{eqnarray*}
0&=&\Suss{1}{u_1}+\cdots+\an{u_1\cdots u_{n-1}}\Suss{1}{u_n}\\
&&+\an{u}\Suss{1}{w_m^{-1}}+\cdots + \an{uw_1^{-1}\cdots w_{3}^{-1}}\Suss{1}{w_2^{-1}}+\cdots +\an{uw_1^{-1}\cdots w_{2}^{-1}}\Suss{1}{w_1^{-1}}\\
&=&\Suss{1}{u_1}+\cdots+\an{u_1\cdots u_{n-1}}\Suss{1}{u_n}\\
&&-\left(    \an{uw_m^{-1}}\Suss{1}{w_m}+\cdots + \an{u w_m^{-1}\cdots w_2^{-1}}\Suss{1}{w_2}+\cdots+\an{u w_m^{-1}\cdots w_1^{-1}}\Suss{1}{w_1}   \right)\\
&=&\Suss{1}{u_1}+\cdots+\an{u_1\cdots u_{n-1}}\Suss{1}{u_n}\\
&&-\left(    \an{w_{1}\cdots w_{m-1}}\Suss{1}{w_m}+\cdots + \an{w_1}\Suss{1}{w_2}+\cdots+\Suss{1}{w_1}   \right)\\
\end{eqnarray*}
as required, since $uw_m^{-1}\cdots w_j^{-1}=w_1\cdots w_{j-1}$ for all $j$.
\item Let $u=\prod_1^n u_i$, $v=\prod_j^m v_j$ in $A^\times$ with $u_i,v_j\in \wn{A}$. Then $uv=u_1\cdots u_n\cdot v_1\cdots v_m$ and hence
\[
\Suss{1}{uv}=\Suss{1}{u}+\an{u}\Suss{1}{v}
\]
by definition.
\end{enumerate}
\end{proof}

\begin{rem}
Observe that property (3) of Corollary \ref{cor:suss1} says that the map $A^\times\to L_3(A)/Z_3(\Sigma_\bullet)$, $u\mapsto \Suss{1}{u}$ is a $1$-cocycle for the group $A^\times$ with values in the $A^\times$-module $L_3(A)/Z_3(\Sigma_\bullet)$.
\end{rem}

\begin{cor} \label{cor:suss2} The conclusions of Corollary \ref{cor:suss1} remain valid if we replace $\Suss{1}{}$ everywhere with $\Suss{2}{}$.
\end{cor}
\begin{proof} If we reverse the roles of the terms $R_3^\alpha(a,b)$ and $R_1^\alpha(a,b)$ in the arguments above, we  derive the same conclusions, but with 
$\Suss{1}{}$ replaced everywhere with $\Suss{2}{}$.
\end{proof}

 If $A$ is a ring in which $\wn{A}$ generates $A^\times$, then by Theorem \ref{thm:l3pres} (1), in  $L_3(A)/Z_3(\Sigma_\bullet)$, we may define  $\cconst{A}:=\cconst{A}(u)$ for any 
$u\in \wn{A}$.

\begin{prop}\label{prop:key}  Let $A$ be a ring which satisfies the condition that $\wn{A}$ generates $A^\times$. Then  
\[
\pf{u}\cconst{A}=\Suss{1}{u}-\Suss{2}{u}\mbox{ in }L_3(A)/Z_3(\Sigma_\bullet)
\]
for all $u\in A^\times$, where $\pf{u}:=\an{u}-1\in \Z[A^\times]$.
\end{prop}

\begin{proof} For any $u\in \wn{A}$, we have
\[
\Suss{1}{u}-\Suss{2}{u}=\an{u}\cconst{A}(u^{-1})-\cconst{A}{u}=\an{u}\cconst{A}-\cconst{A}=\pf{u}\cconst{A} \mbox{ in } L_3(A)/Z_3(\Sigma_\bullet).
\]

On the other hand, for any $u,v\in A^\times$ we have $\pf{uv}=\pf{u}+\an{u}\pf{v}$ and hence, by induction, for $u_1,\ldots,u_n\in A^\times$, 
\[
\pf{u_1\cdots u_n}=\pf{u_1}+\an{u_1}\pf{u_2}+\cdots +\an{u_1\cdots u_{n-1}}\pf{u_n}\mbox{ in } \Z[A^\times].
\]
This, together with Corollaries \ref{cor:suss1} and \ref{cor:suss2}, implies that the required identity extends to all of $A^\times$.
\end{proof}

\subsection{The $\sgr{A}$-module $\hrpb{A}$}
For a  ring $A$, we will introduce for convenience the notation
\[
\hrpb{A}:=\left(\frac{L_3(A)}{Z_3(\Sigma_\bullet)}\right)_{\pspl{2}{A}}\mbox{ and }\hpb{A}=\left(\frac{L_3(A)}{Z_3(\Sigma_\bullet)}\right)_{\pgl{2}{A}}.
\]

The relevance of these modules is the following:

\begin{prop}\label{prop:h3ya} Let $A$ be a ring satisfying
\begin{enumerate}
\item $\wn{A}$ generates $A^\times$, 
\item $Y_4(A)=\emptyset$, and
\item $H_2(|Y(A)|)=0$.
\end{enumerate}
Then there is a natural exact sequence of $\sgr{A}$-modules
\[
H_3(|Y(A)|)_{\pspl{2}{A}}\to\hrpb{A}\to\rpb{A}\to 0,
\]
and a natural sequence of abelian groups
\[
H_3(|Y(A)|)_{\pgl{2}{A}}\to\hpb{A}\to\pb{A}\to 0.
\]
\end{prop}

\begin{proof}  By Theorem \ref{thm:c3}, there is an exact sequence of $\pgl{2}{A}$-modules
\[
H_3(|Y(A)|)\to L_3(A)/Z_2(\Sigma_\bullet)\to Z_2(L_\bullet)\to 0.
\]
Taking $\pspl{2}{A}$-coinvariants and $\pgl{2}{A}$-convariants then gives the required sequences. 
\end{proof}

We consider the condition $Y_4(A):=\emptyset$. In view of the transitivity of $\pgl{2}{A}$ on $2$-simplices, it is enough to consider the existence of $4$-simplices of the form
$\{ \infty,0,1,u,v\}$ with $u,v\in \wn{A}$. It is easy to see that such a set is a $4$-simplex  if and only if $u-v\in A^\times$ and hence if and only if $1-(u/v)\in A^\times$. Thus such $4$-simplices are parametrized by the set $\wn{A}^{(2)}:=\{ \{u,v\}\subset \wn{A}\ |\ u/v\in \wn{A}\}$. Thus $Y_4(A)=\emptyset$, and $L_4(A)=0$ 
if and only if $\wn{A}^{(2)}=\emptyset$. 

\begin{exa} \label{exa:f3} $\wn{\Z/3}=\{ 2=-1\}$ and hence $\wn{\Z/3}^{(2)}=\emptyset$ and thus $L_4(\Z/3)=0$. Since for any ring homomorphism $A\to B$, there is an induced functorial map
$\wn{A}^{(2)}\to\wn{B}^{(2)}$, it follws that $L_4(A)=0$ for any ring $A$ which admits a homomorphism to $\Z/3$.
\end{exa}

\begin{lem}\label{lem:localf3} Let $A$ be a local ring with residue field $\F{3}=\Z/3$. Then $A$ satisfies the three hypotheses of Proposition \ref{prop:h3ya}.
\end{lem}

\begin{proof} We have $A^\times=\{ u\in A^\times\ |\ \bar{u}=1,-1 \mbox{ in }\Z/3\}$ while $\wn{A}=\{ u\in A^\times\ |\ \bar{u}=-1 \mbox{ in }\Z/3\}$. However, if 
$\bar{u}=1$, then $u=(-1)\cdot (-u)$ and $\bar{-1},\bar{-u}\in \wn{A}$. 

By Example \ref{exa:f3}, $\wn{A}^{(2)}=\emptyset$.

Finally, $H_2(|Y(A)|)\cong H_2(L_\bullet(A))$ by Theorem \ref{thm:c3}, and $H_2(L_\bullet(A))=0$ by \cite[Lemma 3.21]{hut:slr}.
\end{proof}

The results above about the $\Z[\pgl{2}{A}]$-module $ L_3(A)/Z_2(\Sigma_\bullet)$ and its special elements imply -- on taking appropriate coinvariants -- corresponding results about the $\sgr{A}$-module
$\hrpb{A}$ and the $\Z$-module $\hpb{A}$. We collect these here:

\begin{lem}\label{lem:hrpbapres} Let $A$ be a ring in which $\wn{A}$ generates $A^\times$. In $\hrpb{A}$, for any $u\in \wn{A}$,  let 
\[
\Suss{1}{u}:=\gpb{u}+\an{-1}\gpb{u^{-1}},\quad \Suss{2}{u}:=\an{u-1}\gpb{u}+\an{1-u^{-1}}\gpb{u^{-1}}.
\]
Also let 
\[
\cconst{A}(u):=\gpb{\frac{1}{1-u}}+\an{-1}\gpb{\frac{1}{1-u^{-1}}}-\Suss{1}{u}.
\]
\begin{enumerate}
\item The $\sgr{A}$-module $\hrpb{A}$ has the following presentation: It is generated by $\gpb{u}$, $u\in \wn{A}$ subject to the relations
\begin{enumerate}
\item $\cconst{A}(u)=\cconst{A}(v)$ for all $u,v\in\wn{A}$,
\item $\an{u}\cconst{A}(u^{-1})-\cconst{A}(u)=\Suss{1}{u}-\Suss{2}{u}$ for all $u\in \wn{A}$, and
\item Whenever $u_1,\ldots,u_n\in \wn{A}$ satisfy $\prod_{i=1}^nu_i=1$ in $A^\times$ then
\[
\Sigma_{i=1}^n\an{u_1\cdots u_{i-1}}\Suss{1}{u_i}=0.
\]
\end{enumerate}
\item The expressions $\Suss{i}{u}$, for $i=1,2$, can be well-defined for all $u\in A^\times$ by
\[
\Suss{i}{u}:=\Sigma_{i=1}^n\an{u_1\cdots u_{i-1}}\Suss{i}{u_i}
\]
whenever $u=\prod_1^n u_i$ with $u_i\in \wn{A}$. They then satisfy
\[
\Suss{i}{uv}=\Suss{i}{u}+\an{u}\Suss{i}{v}
\]
for all $u,v\in A^\times$. 
\item Let $\cconst{A}:=\cconst{A}(u)$ for any $u\in\wn{A}$. The identity
\[
\pf{u}\cconst{A}=\Suss{1}{u}-\Suss{2}{u}
\]
holds for all $u\in A^\times$.
\item For all $u\in A^\times$ and $i\in \{ 1,2\}$, $\an{-1}\Suss{i}{u}=\Suss{i}{u^{-1}}$.
\item $3\cdot\cconst{A}=0$.
\end{enumerate}
\end{lem}

Taking $A^\times$-coinvariants, we therefore also get:

\begin{lem}\label{lem:hpbapres} Let $A$ be a ring in which $\wn{A}$ generates $A^\times$. In $\hpb{A}$,  for any $u\in \wn{A}$,  let 
\[
\Sus{u}:=\gpb{u}+\gpb{u^{-1}}
\]
Also let 
\[
\cconst{A}(u):=\gpb{\frac{1}{1-u}}+\gpb{\frac{1}{1-u^{-1}}}-\Sus{u}.
\]
\begin{enumerate}
\item The $\Z$-module $\hpb{A}$ has the following presentation: It is generated by $\gpb{u}$, $u\in \wn{A}$ subject to the relations
\begin{enumerate}
\item $\cconst{A}(u)=\cconst{A}(v)$ for all $u,v\in\wn{A}$,
\item Whenever $u_1,\ldots,u_n\in \wn{A}$ satisfy $\prod_{i=1}^nu_i=1$ in $A^\times$ then
\[
\Sigma_{i=1}^n\Sus{u_i}=0.
\]
\end{enumerate}
\item The expression $\Sus{u}$ can be well-defined for all $u\in A^\times$ by
\[
\Sus{u}:=\Sigma_{i=1}^n\Sus{u_i}
\]
whenever $u=\prod_1^n u_i$ with $u_i\in \wn{A}$. We have
\[
\Sus{uv}=\Sus{u}+\Sus{v}
\]
for all $u,v\in A^\times$. 
\item For all $u\in A^\times$, $\Sus{u}=\Sus{u^{-1}}$, and hence $\Sus{u^2}=0$.
\item Let $\cconst{A}:=\cconst{A}(u)$ for any $u\in A^\times$. Then  $3\cdot\cconst{A}=0$.
\end{enumerate}
\end{lem}

\subsection{Example: The ring $\F{3}=\Z/3$}
\begin{exa}\label{exa:z3} We consider $A=\Z/3$. Then the space $Y=Y(\Z/3)$ is just a $3$-simplex and $H_i(|Y|)=0$ for all $i>0$. Thus the natural maps $\hrpb{\Z/3}\to\rpb{\Z/3}$ and 
$\hpb{\Z/3}\to \pb{\Z/3}$ are isomorphisms. 

We have $\sgr{\Z/3}$ is the group ring $\Z[\{ 1,\an{-1}\}]$ of a cyclic group of order $2$. $\wn{\Z/3}=\{ -1\}$. So the $\sgr{\Z/3}$-module $\rpb{\Z/3}$ is generated by $\gpb{-1}$. 

We find
$\Suss{1}{-1}=\gpb{-1}+\an{-1}\gpb{-1}=\Suss{2}{-1}$ and $\cconst{\Z/3}=\gpb{-1}+\an{-1}\gpb{-1}-\Suss{1}{-1}=0$. Note that $\an{-1}\Suss{i}{-1}=\Suss{i}{-1}$. 
From $-1\cdot-1=1$, we obtain the relation $\Suss{1}{-1}+\Suss{1}{-1}=2\Suss{1}{-1}=0$; i.e., $2(\gpb{-1}+\an{-1}\gpb{-1})=0$.
 This is the only non-trivial relation satisfied by $\gpb{-1}$. 

Thus $\rpb{\Z/3}\cong \Z[\{ 1,\an{-1}\}]/\an{ 2(1+\an{-1})}$.  As an abelian group $\rpb{\Z/3}\cong \Z\oplus \Z/2$ with generators $\gpb{-1}$ and $\Suss{1}{-1}$.

Taking $(\Z/3)^\times$- coinvariants, it follows that $\pb{\Z/3}$ is cyclic of order $4$ with generator $\gpb{-1}$.

Since $\lambda_1(\gpb{-1})=\pf{-1}^2=-2\pf{-1}$ has infinite order in $\sgr{\Z/3}$ while $\lambda_1(\suss{1}{-1})=0$ it follows that $\rpbker{\Z/3}=\Z/2$ with generator $\Suss{1}{-1}$.
Since the group $B=B(\spl{2}{\Z/3})$ of upper-triangular matrices is nonabelian of order $6$, we have $E^1_{2,0}=H_2(B,\Z)=0$ and hence $E^3_{2,0}=0$, so that $d^3(\Suss{1}{-1})=0$ and hence $\rbl{\Z/3}$ is cyclic of order $2$ with generator $\Suss{1}{-1}$.

Furthermore, $\gpb{-1}\not\in\bl{\Z/3}$ since $\gpb{-1}\not\in \bl{\F{27}}$. Thus $\bl{\Z/3}\subset \pb{\Z/3}$ is cyclic of order $2$ with generator $\Sus{-1}=2\gpb{-1}$, and the natural map $\rbl{\Z/3}\to\bl{\Z/3}$ is an isomorphism.
\end{exa}


\subsection{Example: The ring $\Z/9$}

Since $\Z/9$ is a local ring with residue field $\Z/3=\F{3}$, there is a natural exact sequence
\[
H_3(|Y(\Z/9)|)\to \hrpb{\Z/9}\to\rpb{\Z/9}\to 0.
\]

We begin by determining the structure of $\hrpb{\Z/9}$ as a module over the ring $\sgr{\Z/9}$. Here $\sgr{\Z/9}=\Z[C_2]$ where $C_2$ is cyclic of order $2$ generated by the square class $\an{-1}$.

\begin{prop}\label{prop:hrpbz9}
The $\sgr{\Z/9}$-module $\hrpb{\Z/9}$ is generated by the elements $\gpb{-1}$, $\gpb{2}$ subject to the relations: 
\begin{enumerate}
\item $3(\gpb{2}-\gpb{-1})=0$,
\item $\pf{-1}(\gpb{2}-\gpb{-1})=0$, and
 \item $2\Suss{1}{-1}=2(\gpb{-1}+\an{-1}\gpb{-1})=0$.
\end{enumerate}
Furthermore, $\cconst{\Z/9}=\gpb{2}-\gpb{-1}$ and $\gpb{5}=2\gpb{-1}-\gpb{2}$.
\end{prop}
\begin{proof}
For simplicity, denote the ring $\Z/9$ by $A$. Since $\an{-1}\cconst{A}(-1)=\cconst{A}(-1)$ always (by definition) and since $\cconst{A}(5)=\cconst{A}(2^{-1})=\an{-1}\cconst{A}(2)$
(again, by definition), it follows that the relations (a) in part (1) of Lemma \ref{lem:hrpbapres} reduce  to the single identity $\cconst{A}(2)=\cconst{A}(-1)$.

Since $\an{u}=\an{-1}$ for all $u\in \wn{A}=\{ -1,2,5\}$, we have $\Suss{1}{u}=\Suss{2}{u}$ for $u\in \wn{A}$. Of course, $\pf{-1}\cconst{A}(-1)=0$ also and thus relation
(b) of Lemma \ref{lem:hrpbapres} (1) is redundant.

The relations in $A^\times$ of the elements  $-1,2,5$ are generated by the three relations $-1\cdot -1=1$, $2\cdot 5=1$ and $2^3\cdot -1=1$. By Lemma \ref{lem:hrpbapres}, these give three further relations among the generators $\gpb{-1}$, $\gpb{2}$ and $\gpb{5}$ of $\hrpb{A}$:

From $-1\cdot -1=1$ we get the relation $0=\Suss{1}{-1}+\an{-1}\Suss{1}{-1}=2\Suss{1}{-1}$.

From $2\cdot 5=1$ we derive $0=\Suss{1}{2}+\an{-1}\Suss{1}{5}=2\Suss{1}{2}$.

From $2\cdot 2\cdot 2\cdot -1=1$, we derive $0=\Suss{1}{2}+\an{-1}\Suss{1}{2}+\Suss{1}{2}+\an{-1}\Suss{1}{-1}$. In view of the previous two relations, this gives 
$\Suss{1}{2}=\Suss{1}{-1}$. This relation, in turn, when combined with the relation $2\Suss{1}{-1}=0$ makes the relation $2\Suss{1}{2}=0$ redundant. 

Thus $\hrpb{A}$ is generated by $\gpb{-1}$, $\gpb{2}$ and $\gpb{5}$ subject to the three relations
\begin{enumerate}
\item[(a)] $\cconst{A}(2)=\cconst{A}(-1)$,
\item[(b)] $\Suss{1}{2}=\Suss{1}{-1}$, and
\item[(c)] $2\Suss{1}{-1}=0$.

Then (a) and (b) give $\gpb{-1}+\an{-1}\gpb{2}=\gpb{5}+\an{-1}\gpb{5}=(1+\an{-1})\gpb{5}$. Multiplying both sides by $\pf{-1}$ gives $\pf{-1}(\gpb{-1}-\gpb{2})=0$ -- and hence relation (2) above -- since $\pf{-1}\an{-1}=-\pf{-1}$ and $\pf{-1}(1+\an{-1})=0$.

Then $\cconst{A}=\cconst{A}(2)=\gpb{-1}+\an{-1}\gpb{2}-\suss{1}{2}=\gpb{-1}+\an{-1}\gpb{2}-\suss{1}{-1}=\an{-1}(\gpb{2}-\gpb{-1})=\gpb{2}-\gpb{-1}$ using (b) and (a).
Since $3\cconst{A}=0$, relation (1) follows. Relation (3) always holds. 

Note also that 
\begin{eqnarray*}
\Suss{1}{2}=\Suss{1}{-1}&\iff& \gpb{2}+\an{-1}\gpb{5}=\gpb{-1}+\an{-1}\gpb{-1}\\
&\iff& \an{-1}\gpb{2}+\gpb{5}=\an{-1}\gpb{-1}+\gpb{-1}\\
&\iff & \gpb{5}=\an{-1}(\gpb{-1}-\gpb{2})+\gpb{-1}=(\gpb{-1}-\gpb{2})+\gpb{-1}\mbox{( using (2))}\\
&\iff& \gpb{5}=2\gpb{-1}-\gpb{2}.\\
\end{eqnarray*}

Conversely, given the statement of the Proposition, then  $\cconst{A}(2)=\an{-1}(\gpb{2}-\gpb{-1})=\gpb{2}-\gpb{-1}$ while 
\begin{eqnarray*}
\cconst{A}(-1)&=&\gpb{5}+\an{-1}\gpb{5}-\gpb{-1}-\an{-1}\gpb{-1}=(1+\an{-1})(\gpb{5}-\gpb{-1})\\
&=&(1+\an{-1})(\gpb{-1}-\gpb{2})=-2(\gpb{2}-\gpb{-1})\\
&=&\gpb{2}-\gpb{-1}\mbox{ (using (1).)}\\
\end{eqnarray*}
Thus (a) holds. As noted, (b) holds since $\gpb{5}=2\gpb{-1}-\gpb{2}$. This proves the proposition.
\end{enumerate}

\end{proof}

\begin{cor}\label{cor:hrpbz9} As an additive abelian group, 
\begin{eqnarray*}
\hrpb{\Z/9}&\cong & \Z\oplus \Z/3\oplus \Z/2\mbox{ with generators $\gpb{-1}$, $\cconst{\Z/9}$, $\Suss{1}{-1}$}\\
&\cong &  \Z\oplus \Z/6\mbox{ with generators $\gpb{-1}$, $\bconst{\Z/9}:=\cconst{\Z/9}+\Suss{1}{-1}$}\\
\end{eqnarray*}
\end{cor}

\begin{cor}\label{cor:hrpbz92} Let $k\in \Z/3$. Then for any $u\in \wn{\Z/9}=2+3\Z/3$ we have 
\[
\gpb{(1-3k)u}-\gpb{u}=k\cconst{\Z/9}
\]
in $\hrpb{\Z/9}$.
\end{cor}

\begin{proof} We have $\cconst{\Z/9}=\gpb{2}-\gpb{-1}=\gpb{-1}-\gpb{5}=\gpb{5}-\gpb{2}$. The corollary follows by a direct calculation.
\end{proof}

By Proposition \ref{prop:localhom}, the space $Y(\Z/9)$ is homotopy equivalent to a pointed sum of $16$ $3$-spheres. Furthermore, a basis for $H_3(|Y(\Z/9)|)\cong \Z^{16}$ can be 
described as follows:

For $i\in \{ 0,1,2\}$, let $T_i:=i+3\Z/9$ and let $T_\infty:=\{ \frac{1}{3k}\ |\ k\in \Z/3\}$. Let $x_i=i\in T_i$ for $i\in \{0,1,2\}$ and let $x_\infty:=\infty =\frac{1}{0}\in T_\infty$. 
For $s=(s_0,s_1,s_2,s_\infty)\in S:=\prod_{i\in\projl{\F{3}}}T_i\setminus\{ x_i\}$. Let 
\[
\gamma^s:= \sum_{y\in \prod_i \{ x_i,s_i\}}(-1)^{\delta(y)}(y_0,y_1,y_2,y_\infty)\in L_3(\Z/9).
\]
Then $\{ \gamma^s\}_{s\in S}$ is the required basis.

\begin{lem}\label{lem:rpbz9}  The map $H_3(|Y(\Z/9)|)\to \hrpb{\Z/9}$ is the zero map and hence the map $\hrpb{\Z/9}\to \rpb{\Z/9}$ is an isomorphism.
\end{lem}

\begin{proof} Fix $s=(s_0,s_1,s_2,s_\infty)\in S$ with $s_i\in i+3\Z/3$ and $s_\infty=\frac{1}{3k}$, $k\in \Z/3$.  In $L_3(\Z/9)$, 
\begin{eqnarray*}
\gamma^s&=& \sum (-1)^{\delta(y)}(y_0,y_1,y_2,y_\infty)\\
&=& \sum_{y_\infty=\infty} (-1)^{\delta(y)}(y_0,y_1,y_2,\infty)- \sum_{y_\infty=s_\infty} (-1)^{\delta(y)}(y_0,y_1,y_2,\frac{1}{3k})\\
&=& \sum_{y_\infty=\infty}(-1)^{\delta(y)}\left[ (y_0,y_1,y_2,\infty)-(y_0,y_1,y_2,\frac{1}{3k})\right].\\
\end{eqnarray*}

By Lemma \ref{lem:ya} (2), if $(a,b,c,d)\in L_3(A)$ then $(a,b,c,d)=(\infty,0,1,u)\cdot M_{a,b,c}$ where
\[
u=u(a,b,c,d)=\frac{d(a,c)d(b,d)}{d(b,c)d(a,d)}\in \wn{A}.
\]

It follows that $(a,b,c,d)$ maps to $\an{\det{M_{a,b,c}}}\gpb{u(a,b,c,d)}\in \hrpb{A}$.

(Recall that here  $d(a,b)=a-b$ if $a,b\in A$ and, if $A$ is local,
\[
d\left(\frac{1}{m},a\right)=-d\left( a,\frac{1}{m}\right)=am-1 \mbox{ when } m\in M.)
\]

Now, for $A=\Z/9$, we calculate
\[
u(y_0,y_1,y_2,\infty)=\frac{y_0-y_2}{y_1-y_2} \mbox{ while }u\left( y_0,y_1,y_2,\frac{1}{3k}\right)=(1-3k)\cdot  \frac{y_0-y_2}{y_1-y_2}=(1-3k)\cdot u(y_0,y_1,y_2,\infty)
\]
(using $y_i\cdot 3k=i\cdot 3k \in \Z/9$ when $y_i\in i+3\Z/9$). Thus the image of  $(y_0,y_1,y_2,\infty)-(y_0,y_1,y_2,\frac{1}{3k})$ in $\hrpb{\Z/9}$ is
\begin{eqnarray*}
&&\an{\det{M_{y_0,y_1,y_2}}}\left(\gpb{u(y_0,y_1,y_2,\infty)}-\gpb{(1-3k)\cdot u(y_0,y_1,y_2,\infty)}\right)\\
&=& -\an{\det{M_{y_0,y_1,y_2}}}\cdot k\cconst{\Z/9} \mbox{by Corollary \ref{cor:hrpbz92}}\\
&=&  -k\cconst{\Z/9}\mbox{ since $\an{-1}$ acts trivially on $\cconst{\Z/9}$.}
\end{eqnarray*}
Hence the image of $\gamma^s$ in $\hrpb{\Z/9}$ is 
\begin{eqnarray*}
\left(\sum_{y_\infty=\infty}(-1)^{\delta(y)+1}\right)k\cconst{\Z/9}
\end{eqnarray*}
 and this is $0$ since $\sum_{y_\infty=\infty}(-1)^{\delta(y)+1}=0$.
\end{proof}

\begin{cor} The module $\rbl{\Z/9}$ is cyclic of order $6$ with generator $\bconst{\Z/9}$ and trivial $\sgr{\Z/9}$-action.
\end{cor}

\begin{proof} One has $\lambda_1(\gpb{-1})=\pf{-1}^2=-2\pf{-1}\in \sgr{\Z/9}$ while $\lambda_1(\bconst{\Z/9})=0$. Thus 
$\rpbker{\Z/9}=\Z/6\cdot\bconst{\Z/9}\subset \rpb{\Z/9}$.  Now since $\bconst{\Z[1/2]}\in \rbl{\Z[1/2]}$, it follows that $\bconst{\Z/9}\in\rbl{\Z/9}$. Thus $\rbl{\Z/9}=\rpbker{\Z/9}$.
\end{proof}

\begin{cor} The group $\pb{\Z/9}$ is cyclic of order $12$ with generator $\cconst{\Z/9}+\gpb{-1}$ and $\bl{\Z/9}$ is the subgroup of order $6$ with generator 
$\bconst{\Z/9}=\cconst{\Z/9}+\Sus{-1}$.  In particular, the natural map $\rbl{\Z/9}\to\bl{\Z/9}$ is an isomorphism of groups.
\end{cor}

\begin{proof} The first statement follows from Lemma \ref{lem:rpbz9} together with the fact that $\pb{A}=\rpb{A}_{A^\times/(A^\times)^2}$. 
For the second statement, note that $\gpb{-1}\not\in\bl{\Z/9}$ since $\gpb{-1}\not\in \bl{\F{3}}$. Since $\bconst{\Z/9}\in\bl{\Z/9}$ (it lies in the image of $\rbl{\Z/9}\to\bl{\Z/9}$), the result follows.
\end{proof}
\begin{cor}
If $A$ is any ring admitting a homomorphism to the ring $\Z/9$, then the element $\cconst{A}\in \bl{A}$ has order $3$.
\end{cor}

\begin{rem} As the example of the ring $A=\Z$ shows, it may happen that $\suss{1}{-1}\not \in \rbl{A}$ (although $\suss{1}{-1}\in \rpbker{A}$ always). Recall also, that 
$\cconst{\Z/3}=0$ by Example \ref{exa:z3}.
\end{rem}

\section{A presentation of  $\rpb{A}$ and $\pb{A}$}\label{sec:presrpa}
\subsection{Presentations for $\rpb{A}$ and $\pb{A}$ when $A$ is local  with residue field $\F{3}$}

For a local ring $A$ with residue field $\F{3}$, the standard $2$-parameter five-term dilogarithm type relations are not available in $\rpb{A}$. However, combining Proposition \ref{prop:h3ya} and Lemma \ref{lem:hrpbapres} with the explicit generators of $H_3(|Y(A)|)$ of Proposition \ref{prop:localhom} allows us to write down a presentation of $\rpb{A}$. Apart from the standard relations among $\cconst{A}(u)$, $\Suss{1}{u}$ and $\Suss{2}{u}$, we find a $4$-parameter family of $16$-term relations. (We will show below, by example, that these additional relations are not vacuous.)

\begin{prop}\label{prop:rpbapresf3} Let $A$ be a local ring with maximal ideal $M$ and residue field $\F{3}$. For any $u\in \wn{A}$,  let 
\[
\Suss{1}{u}:=\gpb{u}+\an{-1}\gpb{u^{-1}},\quad \Suss{2}{u}:=\an{u-1}\gpb{u}+\an{1-u^{-1}}\gpb{u^{-1}}.
\]
Also let 
\[
\cconst{A}(u):=\gpb{\frac{1}{1-u}}+\an{-1}\gpb{\frac{1}{1-u^{-1}}}-\Suss{1}{u}.
\]
Then $\rpb{A}$ has the following presentation as a $\sgr{A}$-module:

It is generated by the elements $\gpb{u}$, $u\in \wn{A}$ subject to the relations
\begin{enumerate}
\item $\cconst{A}(u)=\cconst{A}(v)$ for all $u,v\in\wn{A}$,
\item $\an{u}\cconst{A}(u^{-1})-\cconst{A}(u)=\Suss{1}{u}-\Suss{2}{u}$ for all $u\in \wn{A}$, and
\item Whenever $u_1,\ldots,u_n\in \wn{A}$ satisfy $\prod_{i=1}^nu_i=1$ in $A^\times$ then
\[
\Sigma_{i=1}^n\an{u_1\cdots u_{i-1}}\Suss{1}{u_i}=0.
\]
\item Let $(m_0,m_1,m_2,t)\in (M\setminus\{ 0\})^4$ and let $s_i:=i+m_i$ for $i=0,1,2$. Then
\begin{eqnarray*}
0&=&\an{2}\gpb{2}-\an{\frac{1-s_0}{2-s_0}}\gpb{2-s_0}-\an{\frac{2s_1}{2-s_1}}\gpb{\frac{2}{2-s_1}}-\an{\frac{s_2}{s_2-1}}\gpb{\frac{s_2}{s_2-1}}\\
&&-\an{2}\gpb{2(1-t)}+\an{(s_1-s_0)\cdot\frac{s_0-2}{s_1-2}}\gpb{\frac{s_0-2}{s_1-2}}+\an{(1-s_0)\cdot\frac{s_0-s_2}{1-s_2}}\gpb{\frac{s_0-s_2}{1-s_2}}\\
&&+\an{\frac{2-s_0}{1-s_0}}\gpb{2-s_0\cdot\frac{1-t}{1-s_0t}}+\an{s_1\cdot\frac{s_2}{s_2-s_1}}\gpb{\frac{s_2}{s_2-s_1}}+\an{s_1\cdot\frac{2}{2-s_1}}\gpb{(1-s_1t)\cdot
\frac{2}{2-s_1}}\\
&&+\an{\frac{s_2}{s_2-1}}\gpb{(1-t)\cdot \frac{s_2}{s_2-1}}-\an{(s_1-s_0)\cdot\frac{s_0-s_2}{s_1-s_2}}\gpb{\frac{s_0-s_2}{s_1-s_2}}\\
&&-\an{(s_1-s_0)\cdot\frac{s_0-2}{s_1-2}}\gpb{\frac{s_0-2}{s_1-2}\cdot\frac{s_1t-1}{s_0t-1}}-\an{(1-s_0)\cdot\frac{s_0-s_2}{1-s_2}}\gpb{\frac{s_0-s_2}{1-s_2}\cdot\frac{1-t}{1-s_0t}}\\
&&-\an{s_1\cdot\frac{s_2}{s_2-s_1}}\gpb{\frac{s_2}{s_2-s_1}\cdot (1-s_1t)}+\an{(s_1-s_0)\cdot\frac{s_0-s_2}{s_1-s_2}}\gpb{\frac{s_0-s_2}{s_1-s_2}\cdot\frac{s_1t-1}{s_0t-1}}.
\end{eqnarray*}
\end{enumerate}
\end{prop}

\begin{proof} By Proposition \ref{prop:h3ya} we have a natural exact sequence
\[
H_3(|Y(A)|)\to \hrpb{A}\to\rpb{A}\to 0.
\]
 Thus $\rpb{A}$ is generated by the (images of) the elements $\gpb{u}$, $u\in\wn{A}$ subject to the relations defining $\hrpb{A}$ in Lemma \ref{lem:hrpbapres} and the relations given by the images in $\hrpb{A}$ of a set of generators of the abelian group $H_3(|Y(A)|)$. By Proposition \ref{prop:localhom}, a basis for the abelian group $H_3(|Y(A)|)$ is described as follows: Let $T_i:=i+M$ for $i=0,1,2$ and let $T_3:=1/M$. Let $x_i:=i$ for $i=0,1,2$ and let $x_3:=1/0=\infty$. Let $s=(s_0,s_1,s_2,s_3)\in\prod_{i=0^3}(T_i\setminus\{ x_i\})$.
So $s_i=i+m_i$  for $i=0,1,2$ and $s_3=1/t)$ where $(m_0,m_1,m_2,t)\in (M\setminus\{ 0\})^4$. The corresponding basis element of the free abelian group $H_3(Y(A)|)$ is (the class of) 
\[
\gamma^s:=\sum_{y\in \prod_0^3 \{ x_i,s_i\}}(-1)^{\delta(y,s)}y=(0,1,2,\infty)-(s_0,1,2,\infty)-\cdots +(s_0,s_1,s_2,s_3)\in L_3(A)
\]
Recall again that if $(a,b,c,d)\in X_3(A)$, then
$
(a,b,c,d)=(\infty,0,1,u(a,b,c,d))\cdot M_{a,b,c}=\gpb{u}\cdot M_{a,b,c}
$
which maps to 
\[
D(a,b,c,d):=\an{\mathrm{det}(M_{a,b,c})}\gpb{u(a,b,c,d)}=\an{d(a,b)d(b,c)d(c,a)}\gpb{\frac{d(a,c)d(b,d)}{d(a,d)d(b,c)}}
\]
in $\hrpb{A}$. Finally, note that $D(\gamma^s)$ is precisely the expression given in relation (4) above.
\end{proof}

Taking $A^\times$-coinvariants, we deduce a presentation for $\pb{A}$:
\begin{cor}
Let $A$ be a local ring with maximal ideal $M$ and residue field $\F{3}$.  For any $u\in \wn{A}$,  let 
\[
\Sus{u}:=\gpb{u}+\gpb{u^{-1}}
\]
Also let 
\[
\cconst{A}(u):=\gpb{\frac{1}{1-u}}+\gpb{\frac{1}{1-u^{-1}}}-\Sus{u}.
\]
Then $\pb{A}$ has the following presentation as an abelian group:

It is generated by the elements $\gpb{u}$, $u\in \wn{A}$ subject to the relations
\begin{enumerate}
\item $\cconst{A}(u)=\cconst{A}(v)$ for all $u,v\in\wn{A}$,
\item Whenever $u_1,\ldots,u_n\in \wn{A}$ satisfy $\prod_{i=1}^nu_i=1$ in $A^\times$ then
\[
\Sigma_{i=1}^n\Sus{u_i}=0.
\]
\item  Let $(m_0,m_1,m_2,t)\in (M\setminus\{ 0\})^4$ and let $s_i:=i+m_i$ for $i=0,1,2$. Then
\begin{eqnarray*}
0&=&\gpb{2}-\gpb{2-s_0}-\gpb{\frac{2}{2-s_1}}-\gpb{\frac{s_2}{s_2-1}}
-\gpb{2(1-t)}+\gpb{\frac{s_0-2}{s_1-2}}+\gpb{\frac{s_0-s_2}{1-s_2}}\\
&&+\gpb{2-s_0\cdot\frac{1-t}{1-s_0t}}+\gpb{\frac{s_2}{s_2-s_1}}+\gpb{(1-s_1t)\cdot
\frac{2}{2-s_1}}
+\gpb{(1-t)\cdot \frac{s_2}{s_2-1}}\\
&&-\gpb{\frac{s_0-s_2}{s_1-s_2}}
-\gpb{\frac{s_0-2}{s_1-2}\cdot\frac{s_1t-1}{s_0t-1}}-\gpb{\frac{s_0-s_2}{1-s_2}\cdot\frac{1-t}{1-s_0t}}\\
&&-\gpb{\frac{s_2}{s_2-s_1}\cdot (1-s_1t)}+\gpb{\frac{s_0-s_2}{s_1-s_2}\cdot\frac{s_1t-1}{s_0t-1}}.
\end{eqnarray*}
\end{enumerate}
\end{cor}
\subsection{The map $\crpb{A}\to \hrpb{A}$}
\begin{prop}\label{prop:chmap} Let $A$ be a ring with the property that $\wn{A}$ generates $A^\times$. Then there is a natural $\pgl{2}{A}$-module homomorphism $\alpha:Z_2(\Sigma_\bullet)\to L_3(A)/Z_3(\Sigma_\bullet)$ such that the triangle
\[
\xymatrix{
&L_3(A)/Z_3(\Sigma_\bullet)\ar^-{d_3}[d]\\
Z_2(\Sigma_\bullet)\ar[r]\ar^-{\alpha}[ru]&Z_2(L_\bullet)\\
}
\]
of $\pgl{2}{A}$-module homomorphisms  commutes.
\end{prop}

\begin{proof} By Corollary \ref{cor:z2pres}, $Z_2(\Sigma_\bullet)=Z^c\oplus Z^1$ as a $\pgl{2}{A}$-module. 

By Proposition \ref{prop:pres}, $Z^c$ is generated by the element $\cconst{A}$ subject to the relations (1) $\cconst{A}\cdot\omega =\cconst{A}$ and (2) $\cconst{A}\cdot (I+M_{0,1,\infty}+M_{1,\infty,0})=0$. Since $\cconst{A}(u)=\cconst{A}(v)$, for any $u,v\in \wn{A}$,  in $L_3(A)/Z_3(\Sigma_\bullet)$ by Theorem \ref{thm:l3pres} (1), we may define $\alpha(\cconst{A}):=\cconst{A}(u)$ for any $u\in \wn{A}$. Relations (1) and (2) then follow from 
Lemma \ref{lem:cconstids} (2) and (3). 

By  Proposition \ref{prop:pres}, $Z^1$ is generated by the elements $\suss{1}{u}$, $u\in A^\times$, subject to the relations (1) $\suss{1}{u}\cdot \omega =\suss{1}{u^{-1}}$ and (2) 
$\suss{1}{uv}=\suss{1}{u}+\an{u}\suss{1}{v}$ for all $u,v\in A^\times$. For $u\in A^\times$ if $u=\prod_{i=1}^n u_i$ with $u_i\in \wn{A}$. 
We define $\alpha(\suss{1}{u}):=\sum_{i=1}^n\an{u_1\cdots u_{i-1}}\Suss{i}{u_i}$. This is well-defined by Corollary \ref{cor:suss1}. Relation (1) follows from Lemma \ref{lem:cconstids} (1) and relation (2) follows from Corollary \ref{cor:suss1} (3). 

Finally since $d_3(\cconst{A}(u))=\cconst{A}$ and $d_3(\Suss{1}{u})=\suss{1}{u}$ for all $u\in \wn{A}$, the given diagram commutes.
\end{proof}

Taking $\pspl{2}{A}$ and $\pgl{2}{A}$ coinvariants, we deduce:
\begin{cor}\label{cor:chmap}  Let $A$ be a ring with the property that $\wn{A}$ generates $A^\times$.  There is a natural commutative triangle of $\sgr{A}$-modules
\[
\xymatrix{
&\hrpb{A}\ar[d]\\
\crpb{A}\ar[r]\ar^-{\alpha}[ru]&\rpb{A}\\
}
\]
and a commutative triangle of abelian groups
\[
\xymatrix{
&\hpb{A}\ar[d]\\
\cpb{A}\ar[r]\ar^-{\alpha}[ru]&\pb{A}.\\
}
\]
\end{cor}

For convenience, let $\sqcl{A}$ denote the group of square classes $A^\times/(A^\times)^2$ of the ring $A$. Thus $\sgr{A}=\igr{\sqcl{A}}$. Now for $u\in \wn{A}$, we let 
$\sqcl{u}=\sqcl{A,u}$ denote the subgroup of $\sqcl{A}$ generated by the classes $\an{\sigma\cdot u}$, $\sigma\in S_3$.

\begin{lem} Let $A$ be a ring and let $u\in\wn{A}$.
\begin{enumerate}
\item The group $\sqcl{u}$ depends only on the class of $u$ in the orbit set $S_3\backslash \wn{A}$.
\item We always have $\an{-1}\in \sqcl{u}$.
\item The group $\sqcl{u}$ is generated by the three classes $\an{-1}, \an{u}, \an{1-u}$.
\end{enumerate}
\end{lem}

\begin{proof}
\begin{enumerate}
\item This is immediate from the definition.
\item For any $u\in \wn{A}$ we have
\[
-1=u\cdot\left(\frac{u-1}{u}\right)\cdot\frac{1}{1-u}=u\cdot(\rho\cdot u)\cdot(\rho^2\cdot u).
\]
\item For any $v\in\wn{A}$, $\an{\tau\cdot v}=\an{v^{-1}}=\an{v}$. Thus $\sqcl{u}$ is generated by the classes $\an{u}$, $\an{\rho\cdot u}$ and $\an{\rho^2\cdot u}$. However,
$\an{-1}=\an{u}\an{\rho\cdot u}\an{\rho^2\cdot u}$ and hence $\sqcl{u}$ is generated by $\an{-1}$, $\an{u}$ and $\an{\rho^2\cdot u}$. Note that $\an{\rho^2\cdot u}=\an{(1-u)^{-1}}=\an{1-u}$.
\end{enumerate}
\end{proof}

For $u\in \wn{A}$, let $\sgr{u}=\sgr{A,u}$ denote the group ring $\igr{\sqcl{A}/\sqcl{u}}$. Thus there is a natral surjective ring homomorphism $\sgr{A}\to\sgr{u}$ for any $u\in \wn{A}$.

Recall that for any $M\in\pgl{2}{A}$ and for any $u\in\wn{A}$, the image of $\gpb{u}\cdot M$ in $\hrpb{A}$ is $\an{\mathrm{det}(M)}\gpb{u}$. Hence, for any $\sigma\in S_3$, the image of $\sigma\ast\gpb{u}$ in $\hrpb{A}$ is $\an{\mathrm{sgn}(\sigma)}\gpb{\sigma\cdot u}$.

\begin{prop}\label{prop:chrpb} Let $A$ be a ring for which $A^\times$ is generated by $\wn{A}$. Then there is an isomorphism of $\sgr{A}$-modules
\[
\frac{\hrpb{A}}{\crpb{A}}\cong \bigoplus_{\tilde{u}\in S_3\backslash\wn{A}}\sgr{u}\cdot\gpb{\tilde{u}}.
\]
\end{prop}

\begin{proof} Since $\crpb{A}$ is generated by the elements $\suss{1}{u}$, $\cconst{A}$ and $\suss{2}{u}=\suss{1}{u}-\pf{u}\cconst{A}$, it follows from Lemma \ref{lem:hrpbapres} (1), and the definition of the map $\alpha$, that $\coker{\alpha}=\hrpb{A}/\crpb{A}$ is generated by the elements $\gpb{u}$, $u\in\wn{A}$ subject to the relations\\
 (1) $\Suss{1}{u}=(1+\tau)\ast\gpb{u}=0$,\\
 (2) $\cconst{A}(u)=(\rho^2+\tau\rho)\ast\gpb{u}-\Suss{1}{u}=0$, and \\
(3) $\Suss{2}{u}=\an{u-1}\left( \gpb{u}+\an{u}\gpb{u^{-1}}\right)=0$.

Since $\rho^2+\tau\rho=(1+\tau\rho^2)\cdot \rho^2$ and since $\rho^2\ast\gpb{u}=\gpb{\rho^2\cdot u}$, in view of relation (1), relation (2) is equivalent to $
(1+\tau\rho^2)\ast{u}=0$ for all $u\in\wn{A}$.  It follows that $\hrpb{A}/\crpb{A}$ is generated by the elements $\gpb{u}$, subject to the relations $(1+\tau)\ast\gpb{u}=0=
(1+\tau\rho^2)\ast \gpb{u}$ and $\an{u}\gpb{u}+\gpb{u^{-1}}=0$.

For $\sigma\in S_3$ we define $r(\sigma)\in \Z/2$ by $\sgn{\sigma}=(-1)^{r(\sigma)}$. If we set $\epsilon:=-\an{-1}\in \sgr{A}$, the first two identities here are equivalent to $\gpb{u}=\epsilon^{r(\sigma)}\gpb{\sigma\cdot u}$ for $\sigma=\tau$ and $\sigma=\tau\rho^2$. Since these two permutations generate $S_3$, it easily follows that the first two identities are equivalent to the property that 
\[
\gpb{\sigma\cdot u}=\epsilon^{r(\sigma)}\gpb{u} \mbox{ for all }\sigma\in S_3 \mbox{ and }u\in\wn{A}.
\]
The third relation then gives
\[
\an{u}\gpb{u}=-\gpb{u^{-1}}=-\gpb{\tau\cdot u}=-\epsilon\gpb{u}=\an{-1}\gpb{u}
\]
and hence $\an{-u}\gpb{u}=\gpb{u}$ for all $u\in \wn{A}$. For all $\sigma\in S_3$ and $u\in \wn{A}$, this yields more generally the relation
\[
\an{-\sigma\cdot u}\gpb{u}=\epsilon^{r(\sigma)}\an{-\sigma\cdot u}\gpb{\sigma\cdot u}\epsilon^{r(\sigma)}\gpb{\sigma\cdot u}=\epsilon^{2r(\sigma)}\gpb{u}=\gpb{u}.
\]
Since $\an{-1}\in G_u$, this is entirely equivalent to the requirement that $g\gpb{u}=\gpb{u}$ for all $u\in \wn{A}$ and $g\in \sqcl{u}$. 

To conclude: The $\sgr{A}$-module $\hrpb{A}/\crpb{A}$ is generated by the elements $\gpb{u}$, $u\in \wn{A}$ subject to the relations\\
 (a) 
$\gpb{\sigma\cdot u}=\epsilon^{r(\sigma)}\gpb{u}$ for all $u\in\wn{A}$ and for all $\sigma\in S_3$ and\\
 (b) $g\cdot\gpb{u}=\gpb{u}$ for all $u\in\wn{A}$ and for all $g\in \sqcl{u}$. 

It follows that the required isomorphism is induced by the $\sgr{A}$-homomorphism\\
$\beta: \hrpb{A}\to \bigoplus_{\tilde{u}\in S_3\backslash\wn{A}}\sgr{u}\cdot\gpb{\tilde{u}}$ defined by $\beta(\gpb{\sigma\cdot\tilde{u}}):=\epsilon^{r(\sigma)}\gpb{\tilde{u}}$.
\end{proof}
\begin{rem}
Since $\an{-1}=1$ in $\sgr{u}$ for any $u\in\wn{A}$, we have $\epsilon=-1$ in $\sgr{u}$. Hence relation (a) can be simplified to $\gpb{\sigma\cdot u}=\mathrm{sgn}(\sigma)\gpb{u}$ for 
$u\in\wn{A}$ and $\sigma\in S_3$. It follows that the natural map $\beta:\hrpb{A}\to \bigoplus_{\tilde{u}\in S_3\backslash\wn{A}}\sgr{u}\cdot\gpb{\tilde{u}}$ is given by
$\gpb{\sigma\cdot\tilde{u}}\mapsto \mathrm{sgn}(\sigma)\gpb{\tilde{u}}$.
\end{rem}
Since each $\sgr{u}$ is a free abelian group, we deduce:
\begin{cor}
The additive abelian group $\hrpb{A}/\crpb{A}$ is free.
\end{cor}

Taking $A^\times$-coinvariants, we also deduce:
\begin{prop}\label{prop:chpb}  Let $A$ be a ring for which $A^\times$ is generated by $\wn{A}$. There is an isomorphism of abelian groups 
\[
\pb{A}/\cpb{A}\cong \bigoplus_{\tilde{u}\in S_3\backslash\wn{A}}\Z\cdot\gpb{\tilde{u}}
\]
 induced by the map
$\beta:\pb{A}\to \bigoplus_{\tilde{u}\in S_3\backslash\wn{A}}\Z\cdot\gpb{\tilde{u}}$, $\gpb{\sigma\cdot\tilde{u}}\mapsto \mathrm{sgn}(\sigma)\gpb{\tilde{u}}$.
\end{prop}

\begin{exa} By Proposition \ref{prop:chrpb}, Corollary \ref{cor:hrpbz9} and Lemma \ref{lem:rpbz9}, we have a short exact sequence of $\sgr{\Z/9}$-modules
\[
0\to\crpb{\Z/9}\to\hrpb{\Z/9}=\rpb{\Z/9}\to \Z\cdot\gpb{-1}\to 0.
\]
\end{exa}
\section{Some Examples}\label{sec:exa}
\subsection{Example: The ring $\Z[1/2]$}
In \cite[Section 8.4]{hut:bloch1}, the structure of $\rpb{\Z[1/2]}$ was determined by using known information about the structure of $H_3(\spl{2}{\Z[1/2]},\Z)$ as an additive abelian group. In this section, we deduce that the natural homomorphism $\hrpb{\Z[1/2]}\to \rpb{\Z[1/2]}$ is an isomorphism.

\begin{thm}(\cite[Theorem 8.26, Lemma 8.28]{hut:bloch1})\label{thm:rpbzhalf}
As an additive abelian group
\[
\rpb{\Z[1/2]}=\Z/2\oplus\Z/3\oplus\Z/3\oplus\Z\oplus \Z,
\]
where the generators of the factors are $\suss{1}{-1}$, $\cconst{}:=\cconst{\Z[1/2]}$, $\an{2}\cconst{}$, $\gpb{-1}$ and $\gpb{2}$ respectively. Furthermore,
$\an{2}\suss{1}{-1}=\suss{1}{-1}$ and $\cconst{}=\gpb{-1}-\gpb{1/2}$ in $\rpb{\Z[1/2]}$.  Furthermore, $\cconst{}=\gpb{\frac12}-\gpb{-1}$ in $\rpb{\Z[\frac12]}$.
\end{thm}

Now, by  \cite[Proposition 5.26]{hut:f2}, we have $\crpb{\Z[1/2]}=\Z/2\oplus \Z/2\oplus \Z/3\oplus \Z/3\oplus \Z$ with generators $\suss{1}{-1}$, $\an{2}\suss{1}{-1}$, $\cconst{}$, $\an{2}\cconst{}$ and $\suss{1}{2}$ respectively. It follows that the map $\crpb{\Z[1/2]}\to \rpb{\Z[1/2]}$ has kernel of order $2$ generated by $\pf{2}\suss{1}{-1}$. However, we have,

\begin{lem} Let $A$ be a ring in which $2$ is a unit. Then $\pf{2}\suss{1}{-1}=0$ in $\hrpb{A}$.
\end{lem}

This  follows by taking $u=1$  in the more general:
\begin{prop} Let $A$ be a ring and suppose that $u\in A^\times$ satisfies $-u^2\in \wn{A}$. Then $\pf{1+u^2}\suss{1}{-1}=0$ in $\hrpb{A}$.
\end{prop}

\begin{proof} (cf. \cite[Corollary 3.22]{hut:bloch1}) 

For $i=1,2$, we have $\suss{i}{-u^2}=\suss{i}{-1}+\suss{i}{u^2}=\suss{i}{-1}+\pf{u}\suss{i}{-1}=\an{u}\suss{i}{-1}$. Thus $\suss{1}{-u^2}=\an{u}\suss{1}{-1}=\an{u}\suss{2}{-1}=\suss{2}{-u^2}$.

Now
\begin{eqnarray*}
\suss{2}{-u^2}&=& \an{-u^2-1}\gpb{-u^2}+\an{1+u^{-2}}\gpb{-u^{-2}}\\
&=&\an{-(1+u^2)}\left(\gpb{-u^2}+\an{-1}\gpb{-u^{-2}}\right)\\
&=&\an{-(1+u^2)}\suss{1}{-u^2}\\
&=&\an{-(1+u^2}\suss{2}{-u^2}.\\
\end{eqnarray*}

Thus $0=\pf{-(1+u^2)}\suss{2}{-u^2}=\pf{-(1+u^2)}\an{u}\suss{2}{-1}$ and hence\\
 $0=\pf{-(1+u^2)}\suss{2}{-1}=\pf{1+u^2}\suss{i}{-1}$ since $\pf{-1}\suss{i}{-1}=0$.
\end{proof}

\begin{cor} The homorphism $\crpb{\Z[\frac12]}\to\rpb{\Z[\frac12]}$ has kernel of order $2$ with generator $\pf{2}\suss{1}{-1}$.
\end{cor}
\begin{prop} The natural  module homomorphism $\hrpb{\Z[1/2]}\to \rpb{\Z[1/2]}$ is an isomorphsm.
\end{prop}

\begin{proof}
It follows from the calculation above that, if we let ${I}$ denote the image of $\crpb{\Z[1/2]}$ in $\hrpb{\Z[1/2]}$, then the map $I\to \rpb{\Z[1/2]}$ is injective with image, $J$ say, equal to the subgroup generated by 
$\suss{i}{-1}$, $\cconst{}$, $\an{2}\cconst{}$ and $\gpb{2}-\gpb{-1}=\suss{1}{2}-\suss{1}{-1}+\cconst{}$.  

By Proposition \ref{prop:chrpb} and by Theorem \ref{thm:rpbzhalf}, we have an induced isomorphism $\hrpb{\Z[1/2]}/{I}\cong \rpb{\Z[1/2]}/J\cong \Z$ with generator $\gpb{u}$ for any 
$u\in \wn{\Z[1/2]}=\{ -1,2,1/2\}$.
\end{proof}

\begin{rem} By Theorem \ref{thm:c3} there is a complex of $\sgr{\Z[1/2]}$-modules which is exact except possibly at the second term
\[
\xymatrix{H_3(|Y(\Z[1/2])|)_{G}\ar[r] &\hrpb{\Z[1/2]}\ar[r]&\rpb{\Z[1/2]}\ar[r]&H_2(|Y(\Z[1/2])|)_G\ar[r]&0}
\]
where $G:=\pspl{2}{\Z[1/2]}$. It thus follows that $H_2(|Y(\Z[1/2])|)_G=0$ and that the map $H_3(|Y(\Z[1/2])|)\to \hrpb{\Z[1/2]}$ is the zero map. We certainly have 
$H_1(|Y(\Z[1/2])|)=0$ by   \cite[Example 6.13]{hut:ge2arxiv}. It would be interesting to know whether $H_3(|Y(\Z[1/2])|)=0$ and $H_2(|Y(\Z[1/2])|)=0$. More generally, it would be interesting to determine the homotopy type of the $3$-dimensional simplicial complex $Y(\Z[1/2])$.
\end{rem}
\subsection{Example: The ring $\Z_{\an{3}}$}
More generally, we let $A$ be any local ring which contains $\Z$ and has residue field $\F{3}:=\Z/3$. Let $M$ be the maximal ideal of $A$. Thus $M\cap \Z=3\Z$.  Since $|Y(A)|$ is homotopy equivalent to a pointed sum of $3$-spheres, $H_2(|Y(A)|)=0$ and thus there is an exact sequence of $\sgr{A}$-modules 
\[
H_3(|Y(A)|)_{\pspl{2}{A}}\to \hrpb{A}\to\rpb{A}\to 0.
\]
As above, let $D$ denote the map $H_3(|Y(A)|)\to \hrpb{A}$ and for $w\in \wn{A}$ let $\beta_w$ denote the composite
\[
\xymatrix{
\hrpb{A}\ar^-{\beta}[r]&\bigoplus_{u\in S_3\backslash\wn{A}}\sgr{u}\cdot\gpb{u}\ar^-{\pi_w}[r]&\sgr{w}\cdot\gpb{w}\\
}
\]
(where $\pi_w$ is the projection on the $\gpb{w}$-factor).

 Recall (see the proof of Proposition \ref{prop:rpbapresf3}), that $H_3(|Y(A)|)$ has generators of the form $\gamma^s$ where $s=(s_0,s_1,s_2,s_\infty)$ sith $s_i\in i+M$, $i=0,1,2$ and $s_\infty=1/m$ for $m\in M$. For any $n\in \N$, let $s(n):= (3^n,3^n+1,3^n+2,1/3^n)$ and let $\gamma_n:=\gamma^{s(t)}$ be the corresponding element of $H_3(|Y(A)|)$. 

\begin{lem} \label{lem:beta} For $n\geq 1$, let $u_n=u_{n,0}:=2(1-3^n)\in 2+3\Z_{\an{3}}\subset 2+M=\wn{A}$. For all $n,m\geq 2$ we have
\[
\beta_{u_n}(D(\gamma_m))=-\delta_{n,m}\an{2}\gpb{u_n}\in \sgr{u_n}\cdot \gpb{u_n}.
\]
\end{lem}

\begin{proof} As in the proof of Proposition \ref{prop:rpbapresf3} we calculate 
\begin{eqnarray*}
D(\gamma_n)&=&-\an{2}\gpb{2(1-3^n)}+\an{\frac{3^n-2}{3^n-1}}\gpb{\frac{2-3^n}{1+3^n}}+\an{2\cdot\frac{1+3^n}{1-3^n}}\gpb{2\cdot\frac{3^{2n}+3^n-1}{3^n-1}}\\
&&+\an{\frac{3^n+2}{3^n+1}}\gpb{\frac{(2+3^n)(1-3^n)}{1+3^n}}-\an{\frac{3^n-2}{3^n-1}}\gpb{\frac{(3^n-2)(3^{2n}+3^n-1)}{(3^n-1)(3^{2n}-1)}}\\
&&-\an{2\cdot\frac{1-3^n}{1+3^n}}\gpb{\frac{2}{(1+3^n)^2}}-\an{\frac{3^n+1}{3^n+2}}\gpb{(3^n+2)(1-3^n-3^{2n})}\\
&&+\an{2}\gpb{2\cdot\frac{3^{2n}+3^n-1}{3^{2n}-1}}-\an{\frac{3^n-1}{3^n-1}}\left( \gpb{2-3^n}-\gpb{\frac{2-3^n}{1-3^n}}\right)\\
&&-\an{2\cdot\frac{1+3^n}{1-3^n}}\left(\gpb{\frac{2}{1-3^n}}-\gpb{\frac{2}{1+3^n}}\right)+\an{2\cdot\frac{3^n-1}{3^n+1}}\left(\gpb{3^n+2}-\gpb{\frac{3^n+2}{3^n+1}}\right).\\
\end{eqnarray*}

Observe that the expressions within round brackets in the last three terms are all of the form $\gpb{u}-\gpb{u/(u-1)}$, which maps under $\beta$ to $2\gpb{u}$. Thus, applying $\beta$ to 
$D(\gamma_n)$ gives 
\[
\beta(D(\gamma_n))=-\an{2}\gpb{u_{n,0}}+\sum_{i=1}^7\alpha_i\gpb{u_{n,i}}+\sum_{i=8}^{10}2\alpha_i\gpb{u_{n,i}}
\]
where $u_{n,i}\in 2+3\Z_{\an{3}}\subset\wn{A}$ and $\alpha_i=\pm\an{w_i}\in \sgr{A}$.\\
 (Thus $\alpha_1=\an{\frac{3^n-2}{3^n-1}}$ and $u_{n,1}=\frac{2-3^n}{1+3^n}$ etc.)

The $S_3$-orbit of $u_n=2(1-3^n)$ is 
\begin{eqnarray*}
\sigma_1(u_n)=2(1-3^n),\ \sigma_2(u_n)=\frac{1}{2(1-3^n)},\ \sigma_3(u_n)=\frac{2\cdot 3^n-1}{2\cdot 3^n-2}\\
\sigma_4(u_n)=\frac{2\cdot 3^n-2}{2\cdot 3^n-1},\ \sigma_5(u_n)=\frac{1}{2\cdot 3^n-1},\ \sigma_6(u_n)=2\cdot 3^n-1.\\
\end{eqnarray*}

Now, it is entirely elementary to verify, for any $n,m\geq 2$, for any $i=1,\ldots,6$ and for any $j=0,1,\ldots,10$, that 
\[
\sigma_i(u_n)=u_{m,j} \mbox{ if and only if } n=m, \ i=1 \mbox{ and } j=0:
\]

Thus $u_{m,0}=2(1-3^m)$ is a negative integer, and hence the only  term in the $S-3$-orbit it might agree with is $\sigma_1(u_n)=2(1-3^n)$. This clearly forces $n=m$.

Since for $m\geq 2$, $u_{m,1}=\frac{2-3^m}{1+3^m}$ is not an integer, not a unit fraction and is negative it cannot equal $\sigma_i(u_n)$ for any $i$ or $n$. 

Likewise,
 \[
u_{m,2}=2\cdot\frac{3^{2m}+3^m-1}{3^m-1}=\frac{3^{2m}+3^m-1}{(3^m-1)/2}
\]
 is not an integer, not a unit fraction and the numerator and denominator differ by more than $1$. So it cannot be $\sigma_i(u_n)$ for any $n\geq 2$ or $1\leq i\leq 6$.
(However, note that $u_{1,1}=\sigma_2(u_1)=-1/4$.)

Similarly, $u_{m,3}=\frac{(3^m+2)(1-3^m)}{1+3^m}$ and $u_{m,4}=\frac{(3^m-2)(3^{2m}+3^m-1)}{(3^m-1)(3^{2m}-1)}$ are again neither  integers nor a unit fractions and the numerator and denominator differ by more than $1$.

We have $u_{m,5}=\frac{1}{(1+3^m)^2/2}$. This cannot equal $\sigma_5(u_n)=1/(2\cdot 3^n-1)$ since the left-hand denominator is even and the right-hand one is odd.

Likewise, parity considerations rule out  the equality of $u_{m,6}=(3^m+2)(1-3^m-3^{2m})$ and $\sigma_1(u_n)=2(1-3^n)$.

We have $u_{m,7}=\frac{3^{2m}+3^m-1}{(3^m-1)/2}$ which cannot be of the form $\sigma_i(u_n)$ for the same reasons as for $m=2,3,4$.

For all $m\geq 2$, $u_{m,8}=2-3^m$ is a negative integer, not equal to $\sigma_1(u_n)=2(1-3^n)$ for parity reasons. 

For $m\geq 2$, $u_{m,9}=\frac{1}{(1-3^m)/2}$ is a negative unit fraction but cannot equal $\sigma_2(u_n)=\frac{1}{2(1-3^n)}$ since otherwise $3=4\cdot 3^n-3^m$ and here the right-hand side is divisible by $9$. (Note, however, that $u_{2,9}=-1/4=\sigma_2(u_1)$.)

Finally, $u_{m,10}=3^m+2$ is a positive integer, but cannot equal $\sigma_6(u_n)=2\cdot 3^n-1$ since otherwise $3=2\cdot 3^n-3^m$ and the right-hand side again is divisible by $9$.

It now follows that $\beta_{u_n}(\gpb{u_{m,j}})=\delta_{n,m}\delta_{j,0}\gpb{u_n}$ for all $n,m\geq 2$ and for $0\leq j\leq 10$, proving the lemma.
\end{proof}

\begin{prop} Let $A$ be a local ring containing $\Z$ and with residue field $\F{3}$. Let $\gamma_n\in H_3(|Y(A)|)$ as in Lemma \ref{lem:beta}. Then the family of elements
$\{ D(\gamma_n)\}_{n\geq 2}$ in $\ker{\hrpb{A}\to\rpb{A}}$, or in $\ker{\hpb{A}\to\pb{A}}$, are a basis of a free abelian subgroup.
\end{prop}

\begin{proof} This follows immediately from Lemma \ref{lem:beta} since the element $\an{2}\gpb{u_n}$ has infinite order in $\sgr{u_n}\cdot \gpb{u_n}\cong \sgr{u_n}$.
\end{proof}


\end{document}